\documentclass[11pt,twoside]{article}

\usepackage{latexsym}
\usepackage{amsmath}
\usepackage{amsthm}
\usepackage{amssymb}
\usepackage{vmargin}
\usepackage{amscd}
\usepackage{stmaryrd}
\usepackage{euscript}
\usepackage{mathrsfs}
\usepackage{amscd}
\usepackage[all]{xy}
\usepackage{xr}
\DeclareMathAlphabet{\mathpzc}{OT1}{pzc}{m}{it}

\setmargins{32mm}{20mm}{14.6cm}{22cm}{1cm}{1cm}{1cm}{1cm}

\newtheorem{theorem}{Theorem}[section]

\newtheorem{corollary}[theorem]{Corollary}
\newtheorem{lemma}[theorem]{Lemma}
\newtheorem{theorem*}{Theorem}
\newtheorem{proposition*}[theorem*]{Proposition}
\newtheorem{corollary*}[theorem*]{Corollary}
\newtheorem{lemma*}[theorem*]{Lemma}

\theoremstyle{definition}
\newtheorem{definition}[theorem]{Definition}

\newtheorem{definition*}[theorem*]{Definition}

\theoremstyle{remark}
\newtheorem{example}[theorem]{Example}
\newtheorem{examples}[theorem]{Examples}
\newtheorem{remark}[theorem]{Remark}
\newtheorem{remarks}[theorem]{Remarks}

\newtheorem{example*}[theorem*]{Example}
\newtheorem{examples*}[theorem*]{Examples}
\newtheorem{remark*}[theorem*]{Remark}
\newtheorem{remarks*}[theorem*]{Remarks}
\newtheorem{exercise*}[theorem*]{Exercise}

\newcommand\ra{\rightarrow}

\newcommand\id{\mathrm{id}}

\newcommand\ten{\otimes}
\newcommand\vareps{\varepsilon}
\newcommand\eps{\epsilon}

\newcommand\CC{\mathrm{C}}

\renewcommand\H{\mathrm{H}}
\newcommand\z{\mathrm{Z}}

\newcommand\N{\mathbb{N}}
\newcommand\Z{\mathbb{Z}}

\newcommand\vv{\mathbb{V}}

\newcommand\Bu{\mathbb{B}}

\newcommand\bH{\mathbb{H}}

\newcommand\C{\mathcal{C}}

\newcommand\cA{\mathcal{A}}
\newcommand\cB{\mathcal{B}}

\newcommand\cD{\mathcal{D}}
\newcommand\cE{\mathcal{E}}

\newcommand\cL{\mathcal{L}}

\newcommand\cQ{\mathcal{Q}}
\newcommand\cR{\mathcal{R}}

\newcommand\cT{\mathcal{T}}

\newcommand\B{\mathscr{B}}

\newcommand\F{\mathscr{F}}
\newcommand\G{\mathscr{G}}

\renewcommand\O{\mathscr{O}}

\newcommand\sB{\mathscr{B}}
\newcommand\sC{\mathscr{C}}

\newcommand\sE{\mathscr{E}}
\newcommand\sF{\mathscr{F}}
\newcommand\sG{\mathscr{G}}

\newcommand\sM{\mathscr{M}}
\newcommand\sN{\mathscr{N}}
\newcommand\sO{\mathscr{O}}

\newcommand\sT{\mathscr{T}}

\newcommand\sV{\mathscr{V}}

\newcommand\Y{\mathfrak{Y}}

\renewcommand\r{\mathfrak{R}}
\newcommand\Def{\mathfrak{Def}}

\newcommand\fG{\mathfrak{G}}
\newcommand\fH{\mathfrak{H}}

\newcommand\fS{\mathfrak{S}}
\newcommand\fT{\mathfrak{T}}

\newcommand\fX{\mathfrak{X}}
\newcommand\fY{\mathfrak{Y}}
\newcommand\fZ{\mathfrak{Z}}

\renewcommand\L{\Lambda}

\newcommand\m{\mathfrak{m}}

\newcommand\g{\mathfrak{g}}

\newcommand\bu{\ad \Bu}

\renewcommand\hom{\mathscr{H}\!\mathit{om}}

\newcommand\der{\mathscr{D}\!\mathit{er}}
\newcommand\aut{\mathscr{A}\!\mathit{ut}}

\newcommand\Alg{\mathrm{-Alg}}
\newcommand\Mod{\mathrm{-Mod}}
\newcommand\FAlg{\mathrm{-FAlg}}
\newcommand\FMod{\mathrm{-FMod}}
\newcommand\Hom{\mathrm{Hom}}

\newcommand\Ext{\mathrm{Ext}}
\newcommand\EExt{\mathbb{E}\mathrm{xt}}

\newcommand\Der{\mathrm{Der}}

\newcommand\Iso{\mathrm{Iso}}

\newcommand\coker{\mathrm{coker\,}}

\newcommand\Ob{\mathrm{Ob}\,}

\newcommand\Shf{\mathrm{Shf}}

\newcommand\Spec{\mathrm{Spec}\,}

\newcommand\Spf{\mathrm{Spf}\,}

\newcommand\Set{\mathrm{Set}}

\newcommand\Cat{\mathrm{Cat}}

\newcommand\Mor{\mathrm{Mor}\,}

\newcommand\Grp{\mathrm{Grp}}

\newcommand\Grpd{\mathrm{Grpd}}

\newcommand\Cmpts{\mathrm{Cmpts}\,}

\newcommand\ad{\mathrm{ad}}

\newcommand\Lim{\varprojlim}

\newcommand\into{\hookrightarrow}
\newcommand\onto{\twoheadrightarrow}

\newcommand\xra{\xrightarrow}
\newcommand\xla{\xleftarrow}

\newcommand\gd{\mathrm{gr}}

\newcommand\by{\times}
\newcommand\mcl{\mathrm{MC}_L}
\newcommand\mc{\mathrm{MC}}
\newcommand\gl{\mathrm{G}_L}

\newcommand\ddef{\mathrm{Def}}
\newcommand\defl{\mathrm{Def}_L}
\newcommand\defm{\mathrm{Def}_M}

\newcommand\Symm{\mathrm{Symm}}

\newcommand\GL{\mathrm{GL}}

\newcommand\Tot{\mathrm{Tot}\,}

\newcommand\toph{\top_{\mathrm{h}}}
\newcommand\topv{\top_{\mathrm{v}}}
\newcommand\both{\bot_{\mathrm{h}}}
\newcommand\botv{\bot_{\mathrm{v}}}

\renewcommand\iff{\Leftrightarrow}

\newcommand\pd{\partial}

\newcommand\half{\frac{1}{2}}

\newcommand\SDC{S}
\newcommand\Sm{W}

\begin{document}
\title{Deformations via Simplicial Deformation Complexes}
\author{J.P.Pridham\thanks{The author is supported by Trinity College, Cambridge and the Isle of Man Department of Education.}}
\maketitle

%

\section*{Introduction}

There has long been a philosophy that every deformation problem in characteristic $0$ should give rise to a differential graded Lie algebra (DGLA). This DGLA should not only determine the deformation functor (and indeed the deformation groupoid up to equivalence), but also somehow encapsulate more information (corresponding to the higher cohomology) about the original problem. 
We have:
$$
\text{Deformation Problem} \leadsto \text{DGLA} \leadsto \text{Deformation Groupoid} \leadsto \text{Deformation Functor}.
$$

However, there are several respects in which DGLAs are not wholly satisfactory. Finding a DGLA to govern a given problem is not easy  --- Kontsevich describes it as an ``art'' (\cite{Kon}), so the arrow 
$$
\text{Deformation Problem} \leadsto \text{DGLA}
$$
above is merely aspirational. Kontsevich also remarks that, where DGLAs are constructed for geometric problems, the methods used are predominantly analytic. This proves inadequate for the more algebraic problems, such as deformations of an algebraic variety. 
However, in \cite{KS}, Kontsevich and Soibelman go some way towards constructing  DGLAs for algebraic problems, by providing DGLAs to govern deformations over an operad. In non-zero characteristic, DGLAs do not, in general, even give rise to deformation functors.

The reason for these problems is that differential graded objects are  usually natural things to consider only in characteristic zero. This statement is implied by omission in papers such as \cite{QRat}, showing that DG objects suffice when considering homotopy theory over the rationals. In non-zero characteristic, it is necessary to use simplicial objects in their stead. This equivalence of DG and simplicial objects, together with the relative facility of  DG methods, has resulted in a comparative neglect of simplicial methods. As I will demonstrate, there are may scenarios in which simplicial structures arise more naturally. 

In this paper, I define  a simplicial object to replace  DGLAs, namely a Simplicial Deformation Complex (SDC).
This  gives rise to the picture:
$$
\text{Deformation Problem} \leadsto \text{SDC} \leadsto \text{Deformation Groupoid} \leadsto \text{Deformation Functor}.
$$
Here we have the advantage  not only that this works in non-zero characteristic, but also that the first arrow is reasonably canonical. 

Whereas Kontsevich and Soibelman used operads for their notion of algebraicity, I use a weaker notion --- that of monadic adjunctions. Given the existence of a suitable monadic adjunction, an SDC can be constructed to govern the relevant deformation problem. In particular, all operads give rise to monadic adjunctions.

 The frequent existence of such adjunctions, and their  r\^ole in cohomology and homology, are well documented
(e.g.\cite{Mac} Ch.VII \S6, and Notes at the  end of Ch VI).
 It is thus unsurprising that monadic adjunctions should be useful in deformation theory, and they immediately enable the construction of  SDCs for many deformation problems.

One major way in which this approach diverges from previous approaches is that the dual construction is  equally straightforward, providing SDCs whenever there is a suitable comonadic adjunction. As an example, this immediately yields an SDC governing deformations of a local system on an arbitrary topological space.

The real power of this approach, however, lies in the ability to construct SDCs from a combination of  monadic and comonadic adjunctions. For example, to deform  a scheme $X$ is equivalent to deforming its structure sheaf $\sO_X$. The algebra structure of $\sO_X$ can be thought of as monadic, while the sheaf structure is comonadic. This provides an SDC for this example, under more general conditions than those for which a DGLA was constructed  in \cite{Hinich}. By contrast, most previous examples for which DGLAs were constructed were either purely comonadic (e.g. smooth schemes, for which the ring structure does not deform), or purely monadic (e.g. affine schemes, for which the sheaf structure does not deform). Another example of this type is deformation of a Hopf algebra. 

Furthermore, I adapt the results of \cite{QRat} to prove that DGLAs and SDCs are equivalent in characteristic zero. However, while the arrow $\text{DGLA} \leadsto \text{SDC}$ is very natural, the arrow $\text{SDC} \leadsto \text{DGLA}$ is generally not, which helps to explain why constructing DGLAs is an ``art''. 

\tableofcontents

\section{Pre-requisites}
This section  contains a summary of pre-requisite material used throughout the rest of the paper, included for ease of reference.

\subsection{Functors on Artinian Rings}\label{artin}

This section consists merely of a list of well-known results and definitions, included for the convenience of the reader.

Fix $\L$ a complete Noetherian local ring, $\mu$ its maximal ideal, $k$ its residue field. Define  $\widehat{\C_{\L}}$ to be the category of complete local Noetherian  $\L$-algebras with residue field $k$, and $\C_{\L}$ to be the category of Artinian local  $\L$-algebras with residue field $k$

We require that all functors on $\C_{\L}$ satisfy 
\begin{enumerate}
\item[(H0)]$F(k)=\bullet$, the one-point set.
\end{enumerate}
We take the following definitions and results from \cite{Sch}:
%

\begin{definition}\label{xsmallextn} For $p:B \ra A$ in $\C_{\L}$ surjective, $p$ is a small extension if $\ker p=(t)$, a principal ideal, such that $\m_B(t)=(0)$.
\end{definition}

Note that any surjection can be be factorised as a composition of small extensions.

For $F:\C_{\L} \ra \Set$, let $\hat{F}:\widehat{ \C_{\L}} \ra \Set$ by $\hat{F}(R)= \Lim F(R/\m_R^{n})$. 
Note that \mbox{$\hat{F}(R) \xrightarrow{\sim} \Hom(h_R,F)$,} where $h_R:\C_{\L} \ra \Set;\, A \mapsto \Hom(R,A)$.

\begin{definition}
We will say a functor $F:\C_{\L} \ra \Set$ is pro-representable if it is isomorphic to $h_R$, for some $R \in \widehat{\C_{\L}}$. By the above remark, this isomorphism is determined by an element $\xi \in \hat{F}(R)$. We say the pro-couple $(R,\xi)$ pro-represents $F$.
\end{definition}

\begin{definition} A natural transformation $\phi:F \ra G$ in $[\C_{\L},\Set]$ is called:
\begin{enumerate}
\item unramified if $\phi:t_F \ra t_G$ is injective, where $t_F=F(k[\epsilon])$.
\item smooth if for every $B \twoheadrightarrow A$ in $\C_{\L}$, we have $F(B) \twoheadrightarrow G(B)\by_{G(A)}F(A)$
\item \'etale if it is smooth and unramified.
\end{enumerate}
\end{definition}

\begin{definition}
$F:\C_{\L} \ra \Set$ is smooth if and only if $F \ra \bullet$ is smooth.
\end{definition}

\begin{definition} A pro-couple $(R,\xi)$ is a hull for $F$ if the induced map $h_R \ra F$ is \'etale.
\end{definition}

\begin{theorem} Let $(R, \xi)$, $(R',\xi')$ be hull of $F$. Then there exists an isomorphism $u:R \ra R'$ such that $F(u)(\xi)=\xi'$.
\end{theorem}

\begin{lemma} Suppose $F$ is a functor such that $$F(k[V]\by_k k[W]) \xrightarrow{\sim}F(k[V])\by F(k[W])$$ for vector spaces $V$ and $W$ over $k$, where $k[V]:=k\oplus V \in \C_{\L}$ in which \mbox{$V^{2}=0$.} Then $F(k[V])$, and in particular $t_F$, has a canonical vector space structure, and \mbox{$F(k[V]) \cong t_F \ten V$.}
\end{lemma}

\begin{theorem}\label{Sch}
Let $F:\C_{\L} \ra \Set$. Let $A' \ra A$ and $A'' \ra A$ be morphisms in $\C_{\L}$, and consider the map:
\begin{enumerate}
\item[$(\dagger)$]  $F(A'\by_A A'') \ra F(A')\by_{F(A)}F(A'').$
\end{enumerate}
Then
\begin{enumerate}
\item $F$ has a hull if and only if $F$ has properties (H1), (H2) and (H3) below:
\begin{enumerate}
\item[{\rm (H1)}] $(\dagger)$ is a surjection whenever $A'' \ra A$ is a small extension.
\item[{\rm(H2)}] $(\dagger)$ is a bijection when $A=k,\quad A''=k[\epsilon]$.
\item[{\rm(H3)}] $\dim_k(t_F) < \infty$.
\end{enumerate}
\item $F$ is pro-representable if and only if $F$ has the additional property (H4):
\begin{enumerate}
\item[{\rm(H4)}]
 $F(A'\by_A A'') \xrightarrow{\sim} F(A')\by_{F(A)}F(A''). $
\end{enumerate}
for any small extension $A' \ra A$.
\end{enumerate}
\begin{proof}
\cite{Sch}, Theorem 2.11.
\end{proof}
\end{theorem}

\begin{definition}
$F:\C_{\L} \ra \Set$ is homogeneous if $\eta : F(B\by_AC)\ra F(B)\by_{F(A)}F(C)$ is an isomorphism for every $B \twoheadrightarrow A$.

Note that a homogeneous functor satisfies conditions (H1), (H2) and (H4).
\end{definition}

\begin{definition}
$F:\C_{\L} \ra \Set$ is a deformation functor if:
\begin{enumerate}
\item $\eta$ is surjective whenever $B \twoheadrightarrow A$.
\item $\eta$ is an isomorphism whenever $A=k$.
\end{enumerate}
Note that a deformation functor satisfies conditions (H1) and (H2).
\end{definition}

The following results are proved by Manetti (in \cite{Man}):

\begin{theorem}\label{SSC} (Standard Smoothness Criterion) Given $\phi :F \ra G$ a morphism of deformation functors, with \mbox{$(V,v_e) \xrightarrow{\phi'}(W,w_e)$} a compatible morphism of obstruction theories, if $(V,v_e)$ is complete, $V \xrightarrow{\phi'} W$ injective, and $t_F \ra t_G$ surjective, then $\phi$ is smooth.
\begin{proof} \cite{Man}, Proposition 2.17.
\end{proof}
\end{theorem}

For functors $F:\C_{\L} \ra \Set$ and $G:\C_{\L} \ra \Grp$, we say that $G$ acts on $F$ if we have a functorial group action $G(A) \by F(A) \xra{*} F(A)$, for each $A$ in $\C_{\L}$.  The quotient functor $F/G$ is defined by $(F/G)(A)=F(A)/G(A)$.

\begin{theorem}\label{Man1} If $F:\C_{\L} \ra \Set$, a deformation functor, and $G:\C_{\L} \ra \Grp$ a smooth deformation functor, with $G$ acting on $F$, then $D:=F/G$ is a
deformation functor, and if $\nu :t_G \ra t_F$ denotes $h \mapsto h*0$, then $t_D=\coker\nu$, and the universal obstruction theories of $D$ and $F$ are isomorphic.
\begin{proof} \cite{Man}, Lemma 2.20.
\end{proof}
\end{theorem}

\begin{theorem}\label{Man2} For $F:\C_{\L} \ra \Set$ is homogeneous, and $G:\C_{\L} \ra \Grp$ a deformation functor, given $a,b \in F(R)$, define $\Iso(a,b): \C_R \ra \Set$  by 
$$ 
\Iso(a,b)(R \xrightarrow{f}A)=\{g \in G(A) | g*f(a)=f(b)\}.
$$ 
Then $\Iso(a,b)$ is a deformation functor, with tangent space $\ker\nu$ and, if $G$ is moreover smooth, complete obstruction space $\coker\nu=t_D$. Furthermore, if $G$ is homogeneous, then so is $\Iso(a,b)$.
\begin{proof} \cite{Man}, Proposition 2.21.
\end{proof}
\end{theorem}

\begin{theorem}\label{Man3} If $G,G'$ smooth deformation functors, acting on $F,F'$ respectively, with $F,F'$ homogeneous, $\ker \nu \ra \ker \nu'$ surjective, and $\coker\nu\ra\coker\nu'$ injective, then $F/G \ra F'/G'$ is injective.
\begin{proof} \cite{Man}, Corollary 2.22.
\end{proof}
\end{theorem}

This final result does not appear in \cite{Man}, but proves extremely useful:

\begin{corollary}\label{keyhgs}
If $F:\C_{\L} \ra \Set$ and $G:\C_{\L} \ra \Grp$ are deformation functors, with $G$ acting on $F$, let $D:=F/G$, then:
\begin{enumerate}
\item If $G$ is smooth, then $\eta_D$ is surjective for every $B \onto A$ (i.e. $D$ is a deformation functor).
\item If $F$ is homogeneous and $\ker \nu=0$, then $\eta_D$ is injective for every $B \onto A$.
\end{enumerate}
Thus, in particular, $F/G$ will be homogeneous if $F$ is homogeneous, $G$ is a smooth deformation functor and $\ker \nu=0$.\begin{proof}
The first part is just Theorem \ref{Man1}. For the second part, assume we are given 
$$
\alpha, \beta \in D(B\by_A C) \text{  such that } \eta(\alpha)=\eta(\beta),
$$
so $\alpha|_B=\beta|_B,\quad \alpha|_C=\beta|_C.$ 
Lift $\alpha,\beta$ to $\tilde{\alpha},\tilde{\beta} \in F(B\by_A C)$. We have \mbox{$g_B \in G(B),\quad g_C \in G(C)$} such that 
$$
g_B*\tilde{\alpha}|_B=\tilde{\beta}|_B,\quad g_C*\tilde{\alpha}|_C=\tilde{\beta}|_C.
$$ 
Thus 
$
g_B|_A, g_C|_A \in \Iso(\tilde{\alpha}|_A,\tilde{\beta}|_A)(A)$. By Theorem \ref{Man2}, $\Iso(\tilde{\alpha}|_A,\tilde{\beta}|_A)$ has a hull, with tangent space $\ker \nu=0$, so $\Iso(\tilde{\alpha}|_A,\tilde{\beta}|_A)(A)$ contains at most one element, and \mbox{$g_B|A=g_C|A$.} Hence 
$$
(g_B,g_C) \in G(B)\by_{G(A)}G(C),
$$
 and we may use surjectivity of $\eta_G$ to lift this to $g \in G(B\by_A C)$. Now $\eta_F(g*\tilde{\alpha})=\eta_F(\tilde{\beta})$, so homogeneity of $F$ implies that $g*\tilde{\alpha}=\tilde{\beta}$, so $\alpha=\beta$.
\end{proof}
\end{corollary}

To summarise the results concerning the pro-representability of the quotient \mbox{$D=F/G$,} we have:
\begin{enumerate}
\item If $F$ is a deformation functor and $G$ a smooth deformation functor, with $\coker \nu$ finite dimensional, then $D$ has a hull.
\item If $F$ is homogeneous and $G$ a smooth deformation functor, with $\ker \nu =0$ and $\coker \nu$ finite dimensional, then $D$ is pro-representable.
\end{enumerate}

\subsection{Differential Graded Lie Algebras}\label{xdgla}
The results of this section can all be found in \cite{Man}.
\begin{definition} A DGLA over a field $k$  of characteristic $0$ is a  graded vector space $L=\bigoplus_{i \in \Z} L^i$ over $k$, equipped with operators $[,]:L \by L \ra L$ bilinear and $d:L \ra L$ linear,  satisfying:

\begin{enumerate}
\item $[L^i,L^j] \subset L^{i+j}.$

\item $[a,b]+(-1)^{\bar{a}\bar{b}}[b,a]=0$.

\item $(-1)^{\bar{c}\bar{a}}[a,[b,c]]+ (-1)^{\bar{a}\bar{b}}[b,[c,a]]+ (-1)^{\bar{b}\bar{c}}[c,[a,b]]=0$.

\item $d(L^i) \subset L^{i+1}$.

\item $d \circ d =0.$

\item $d[a,b] = [da,b] +(-1)^{\bar{a}}[a,db]$
\end{enumerate}

Here $\bar{a}$ denotes the degree of $a$, mod $ 2$, for $a$ homogeneous.

\end{definition}

\begin{remark}
 Note that this definition of a DGLA does allow for negative degrees; however, it is natural to require that $\H^{<0}(L\ten k)=0$, and in Section \ref{sdcdgla} it will be shown that all the deformation problems we consider can be described by DGLAs with $L^{<0}=0$.
\end{remark}

Fix a  DGLA $L$.

\begin{definition} The Maurer-Cartan functor 
$\mcl:\C_{k} \ra \Set$ is defined by 
$$ \mcl(A)=\{x \in L^1\otimes \m_A |dx+\half[x,x]=0\}.$$
\end{definition}

Observe that for $\omega \in L^1\ten \m_A$, 
$$
d\omega +\half[\omega,\omega]=0 \Rightarrow (d+\ad_{\omega})\circ(d+\ad_{\omega})=0,
$$ 
so $(L\ten A, [,], d+\ad_{\omega})$ is a DGLA over $A$.

\begin{definition}
Define the gauge functor $\gl:\C_k \ra \Grp$ by $\gl(A)=\exp(L^0 \ten \m_A)$. 
\end{definition}

We may define another DGLA, $L_d$, by 
$$
L_d^i=\left\{\begin{matrix} L^1\oplus k d	& i=1\\
                           L^i 			& i \ne 1, \end{matrix}\right.
$$
 with 
$$
d_d(d)=0,\quad [d,d]=0,\quad [d,a]_d=da,\quad \forall a \in L.
$$ 

\begin{lemma} 
 $exp(L^0 \ten \m_A)$ commutes with $[,]$ when acting  on $L_d \ten A$ via the adjoint action.
\end{lemma}

\begin{corollary}
Since $\exp(L^0 \ten \m_A)$ preserves $L^1\ten \m_A+d \subset L_d\ten A$ under the adjoint action, and 
$$
x \in \mcl(A) \iff [x+d,x+d]=0,
$$
 the adjoint action of $\exp(L^0 \ten \m_A)$ on $L^1\ten A+d$ induces an action of $\gl(A)$ on $\mcl(A)$, which we will call the gauge action.
\end{corollary}

\begin{definition} $\defl=\mcl/\gl$, the quotient being given by the gauge action \mbox{$\alpha(x)= \ad_{\alpha}(x+d)-d$.}
  Observe that $\gl$ is homogeneous. 

 For $a \in \mcl(A)$, define $K_a:\C_{A} \to \Grp$ by 
$$
K_a(B)=\exp((d+\ad_a)L^{-1}\ten \m_B).
$$
 Note that this makes sense, since $(d+\ad_a)^2=0$, so 
$$
(d+\ad_a)L^{-1}\ten \m_B \le L^0\ten \m_B
$$
is a Lie subalgebra.
Note that $K_a$ is then a subfunctor of $\Iso(a,a)\le \gl(A)$, so acts on $\Iso(a,b)$ by right multiplication.

Define the deformation groupoid $\mathfrak{Def}_L$ to have objects $\mcl$, and morphisms  given by $\Iso(a,b)/K_a$.
\end{definition}

Now, 
$
t_{\gl}=L^0\epsilon ,
$ 
and 
$
t_{\mcl}=\z^1(L \epsilon),
$
 with action 
\begin{eqnarray*}
t_{\gl}\by t_{\mcl} &\ra& t_{\mcl};\\ 
(b,x) &\mapsto& x+db, \text{ so }
\end{eqnarray*}
$$
 t_{\defl}=\H^1(L). 
$$  

\begin{lemma} $\H^2(L)$ is a complete obstruction space for
$\mcl$.
\begin{proof} Given a small extension 
$$
e:0 \ra J \ra A \ra B\ra 0,
$$ and $x \in \mcl(B)$, lift $x$ to $\tilde{x} \in L^1\ten \m_A$, and let 
$$
h=d\tilde{x} +\half [\tilde{x},\tilde{x}] \in L^2\ten \m_A.
$$
 In fact, $h \in L^2\ten J$, as $dx+\half[x,x]=0$. 

Now,
$$
dh=d^2\tilde{x}+[d\tilde{x},\tilde{x}]= [h-\half [\tilde{x},\tilde{x}],\tilde{x}]=[h,\tilde{x}]=0,
$$ 
since $[[\tilde{x},\tilde{x}],\tilde{x}]=0$ and $J.\m_A=0$. Let 
$$
v_e(x)=[h]\in \H^2(L\ten J)=\H^2(L)\ten_k J.
$$
This is well-defined: if $y=\tilde{x}+z$, for $z \in L^1\ten J$, then
$$
dy +\half [y,y]=d\tilde{x}+dz+\half[\tilde{x},\tilde{x}]+\half[z,z]+[\tilde{x},z]=h+dz,
$$ 
as $J.\m_A=0$.

This construction is clearly functorial, so it follows that $(\H^2(L),v_e)$ is a complete obstruction theory
for $\mcl$.
\end{proof}
\end{lemma}

Now Theorem \ref{Man1} implies the following:

\begin{theorem} $\defl$ is a deformation functor, $t_{\defl} \cong \H^1(L)$, and $\H^2(L)$ is a complete obstruction theory for $\defl$.
\end{theorem}

The other theorems can be used to prove:

\begin{theorem}\label{qis}
If $\phi :L \ra M$ is a morphism of DGLAs, and 
$$\H^i(\phi):\H^i(L) \ra \H^i(M)$$ are the induced maps on cohomology, then:
\begin{enumerate}
\item If $\H^1(\phi)$ is bijective, and $\H^2(\phi)$ injective, then $\defl \ra \defm$ is \'etale.
\item If also $\H^0(\phi)$ is surjective, then  $\defl \ra \defm$ is an isomorphism.
\item Provided condition 1 holds, $\mathfrak{Def}_L \to \mathfrak{Def}_M$ is an equivalence of functors of  groupoids if and only if $\H^0(\phi)$ is an isomorphism.
\end{enumerate}
\begin{proof}
\cite{Man}, Theorem 3.1, mutatis mutandis.
\end{proof}
\end{theorem}

\begin{theorem}
If $\H^0(L)=0$, then $\defl$ is homogeneous.
\begin{proof}
This is essentially Theorem \ref{keyhgs}, with some allowance made for $K_a$.
\end{proof}
\end{theorem}

Thus, in particular, a quasi-isomorphism of DGLAs gives an isomorphism of deformation functors and deformation groupoids.

\begin{example} Deformation of algebras (due to \cite{Kon}).

We wish to  describe flat deformations $R_A/ A$ of any  algebra $R_0/k$:
$$
\begin{CD}
X=\Spec R_0 @>>> \Spec R_A= X_A\\
 @V{\text{flat}}VV  @VV{\text{flat}}V  \\
  \Spec k @>>> \Spec A.  
\end{CD}
$$
Explicitly, the deformation groupoid has objects flat $R_A/ A$, such that \mbox{$R_A \ten_A k =R_0$,}  with morphisms given by infinitesimal isomorphisms of $R_A$, i.e. $ A$-algebra isomorphisms $\theta$ such that \mbox{$\theta \equiv \id \mod \m_A$.}

Take any free resolution $R_{\bullet} \to R$ of $R$, i.e. a free graded algebra $R_*$,  with a differential $\delta$ of degree $1$, such that 
$$
\H_i(R_{\bullet})= \left\{ \begin{matrix} R &i=0 \\  0 &i > 0. \end{matrix} \right.
$$

Now, let 
$$
L^*=\Der_k(R_*,R_*),
$$
with differential $\ad_{\delta}$.

Let $\Def_{R_0/k}(A)$ denote the groupoid of deformations  $R_A/A$ of $R_0$, with morphisms given by infinitesimal isomorphisms. Given $\omega \in \mcl(A)$, it follows that $(\delta+\omega)^2=0$. Define a map
$$
\Def_L(A) \xra{\Phi} \Def_{R_0/k}(A),
$$
sending $\omega \in \mcl(A)$ to the flat $A$-algebra 
$$
H_0(R_*\ten A, \delta+\omega).
$$
 In fact, $\Phi$ gives an essential surjection.

Next, observe that $\gl(A)$ consists precisely of infinitesimal $A$-algebra automorphisms of $R^*\ten A$, and that the gauge action of $\theta$ sends $\delta+\omega$ to $\theta \circ (\delta+\omega) \circ \theta^{-1}$.

Now,  $\theta \in \Hom(\omega, \omega')$ maps under $\Phi$ to
$$
H_0(\theta):H_0(R_*\ten A, \delta+\omega) \to H_0(R_*\ten A, \delta+\omega'),
$$
and it turns out that $\Phi$ is full.

To see that $\Phi$ is  faithful, note that for  $\theta \in \Iso(\omega,\omega)(A)$, 
$$
H_0(\theta):H_0(R_*\ten A, \delta+\omega) \to H_0(R_*\ten A, \delta+\omega)
$$
is the identity if and only if $\theta \in K_{\omega}(A)$.
\end{example}

\subsection{Monadic Adjunctions}\label{xmon}

This section just recalls some standard definitions concerning adjunctions, and fixes notation which will be used throughout the paper.

\begin{definition}
For categories $\cD, \cE$, and a pair of functors
$$
\xymatrix@1{\cD \ar@<1ex>[r]^G & \cE \ar@<1ex>[l]^F},
$$
recall that an adjunction $F \dashv G$ is a natural isomorphism
$$
\Hom_{\cD} (FA,B) \cong \Hom_{\cE}(A,GB).
$$
We say that $F$ is left adjoint to $G$, or $G$ is right adjoint to $F$. Let $\bot = FG$, and $\top=GF$. To give an adjunction is equivalent to giving two natural transformations, the unit
$$
\eta:\id_{\cE} \to  \top,
$$
and the co-unit
$$
\vareps: \bot \to \id_{\cD},       
$$
satisfying the triangle identities:
$$
\xymatrix{F \ar[r]^{F\eta} \ar[dr]_{\id_F} &FGF \ar[d]^{\vareps F} \\
					&F}			
\text{ and }
\xymatrix{G \ar[r]^{\eta G} \ar[dr]_{\id_G} &GFG \ar[d]^{G\vareps } \\
					&G.}
$$
In terms of the unit and co-unit, the adjunction is given as follows:
\begin{eqnarray*}			
\Hom_{\cD} (FA,B) &\cong& \Hom_{\cE}(A,GB)\\
\theta		&\mapsto& G\theta \circ \eta_A\\
\vareps_B \circ F\phi &\mapsfrom& \phi.
\end{eqnarray*}
Note that the $\eta_A$ and $\vareps_B$ can be recovered from the adjunction, since they correspond under the adjunction to the identity maps on $FA$ and $GB$, respectively. 
\end{definition}
	
\begin{examples}
\begin{enumerate}
\item
Given a morphism $X \xra{f} Y$ of schemes, there is an adjunction 
$$
\xymatrix@1{\Shf(X) \ar@<1ex>[r]^{f_*} & \Shf(Y) \ar@<1ex>[l]^{f^{-1}}_{\bot}},
$$
between the category of sheaves on $Y$ and the category of sheaves on $X$, i.e. a natural isomorphism
$$
\Hom_X(f^{-1}\sG,\sF) \cong \Hom_Y(\sG, f_*\sF),
$$
determined by the canonical maps
$$
\vareps_{\sF}:f^{-1}f_*\sF \to \sF, \quad  \eta_{\sG}: \sG \to f_*f^{-1}\sG.
$$
\item
There is also an adjunction
$$
\xymatrix@1{\O_X\Mod \ar@<1ex>[r]^{f_*} & \O_Y\Mod \ar@<1ex>[l]^{f^*}_{\bot}},
$$
$$
\Hom_{\O_X}(f^{*}\sG,\sF) \cong \Hom_{\O_Y}(\sG, f_*\sF),
$$
$$
\vareps_{\sF}:f^{*}f_*\sF \to \sF,\quad \eta_{\sG}: \sG \to f_*f^{*}\sG.
$$
\end{enumerate}
\end{examples}

If  a functor has a left adjoint, then it preserves all (inverse) limits. Conversely, provided the categories involved satisfy various additional conditions, the Special Adjoint Functor Theorem (\cite{Mac} Ch. V.8) proves that any functor which preserves inverse limits has a left adjoint.

Given an adjunction
$$
\xymatrix@1{\cD \ar@<1ex>[r]^U_{\top} & \cE \ar@<1ex>[l]^F}
$$
with unit $\eta:\id \to UF$ and co-unit $\vareps:FU \to \id$, we let $\top=UF$, and define the category of $\top$-algebras, $\cE^{\top}$, to have objects
$$
\top E \xra{\theta} E,
$$
such that the following diagrams commute:
$$
\xymatrix{ E  \ar[r]^{\eta_E} \ar[dr]_{\id}   &\top E \ar[d]^{\theta} \\
				&  E }\quad \text{ and }\quad\quad
\xymatrix{ \top^2 E \ar[r]^{\top \theta} \ar[d]_{U\vareps_{FE}}   &\top E \ar[d]^{\theta}  \\
		\top E \ar[r]^{\theta}  &E .}
$$	
A morphism 
$$
g: \left( \begin{CD} \top E_1 \\ @VV{\theta}V \\ E_1 \end{CD} \right) \to  \left( \begin{CD} \top E_2 \\ @VV{\phi}V \\ E_2 \end{CD} \right) 
$$	
of $\top$-algebras is a morphism $g:E_1 \to E_2$ such that
$$
\begin{CD} \top E_1 @>{\top g}>>  \top E_2 \\
		@V{\theta}VV  @VV{\phi}V \\
		E_1 @>>{g}>  E_2 \end{CD}
$$
commutes.

We define the comparison functor $K:\cD \to \cE^{\top}$ by
$$
B \mapsto \left(\begin{CD} UFUB \\ @VV{U\vareps_B}V \\ UB \end{CD}\right)
$$
on objects, and $K(g)=U(g)$ on morphisms.

\begin{definition}
The adjunction 
$$
\xymatrix@1{\cD \ar@<1ex>[r]^U_{\top} & \cE \ar@<1ex>[l]^F},
$$
is said to be  \emph{monadic} if $K:\cD \to \cE^{\top}$ is an equivalence.
\end{definition}

Intuitively, monadic adjunctions correspond to algebraic theories, such as the adjunction
$$
\xymatrix@1{k\Alg \ar@<1ex>[r]^U_{\top} &k\mathrm{-Vect}  \ar@<1ex>[l]^{\Symm_k}},
$$
 between $k$-algebras and vector spaces over $k$, $U$ being the forgetful functor. Other examples are Lie algebras  or non-commutative algebras, with their respective free functors, over vector spaces, or vector spaces over sets.

\begin{example} Probably the clearest example of a monadic adjunction is that between sets and $G$-sets. Given a group $G$, The category $G\mathrm{-Set}$ consists of sets $S$, with a $G$-action $G \by S \to S$. We have the adjunction:
$$
\xymatrix@1{G\mathrm{-Set} \ar@<1ex>[r]^-U_-{\top} &\mathrm{Set}  \ar@<1ex>[l]^-{G\by-}},
$$
with $\eta_S = (e,\id):S \to G \by S$, and $\vareps_T = *:G \by T \to T$.

A $\top$-algebra consists of a set $S$, together with a map $\top S = G\by S \xra{f} S$, such that:
$$
\xymatrix{ S  \ar[r]^-{(e,\id_S)} \ar[dr]_{\id_S}   &G \by S \ar[d]^{f} \\
				&  S }\quad \text{ and }\quad\quad
\xymatrix{ G \by G \by S \ar[r]^-{(\id_G,f)} \ar[d]_{(*,\id_S)}   &G \by S \ar[d]^{f}  \\
		G \by S \ar[r]^{f}  &S }
$$	
commute, which is precisely the same as a $G$-set, hence the equivalence.
\end{example}

Dually, if we have an adjunction 
$$
\xymatrix@1{\cD \ar@<1ex>[r]^V & \cE \ar@<1ex>[l]^G_{\bot}},
$$
with co-unit $\gamma:VG \to \id$ and  unit $\alpha:\id \to GV$, we let $\bot=VG$, and can define the category of $\bot$-coalgebras, $\cE^{\bot}$, with a functor $K:\cD \to \cE^{\bot}$.

\begin{definition}
The adjunction
$$
\xymatrix@1{\cD \ar@<1ex>[r]^V & \cE \ar@<1ex>[l]^G_{\bot}},
$$
is said to be comonadic if $K:\cD \to \cE^{\bot}$ is an equivalence.
\end{definition}

Finally, the equivalence version of Beck's monadicity theorem gives a criterion for an adjunction to be monadic:

\begin{definition}\label{split}
A split fork is a diagram 
$$
\xymatrix@1{A \ar@<0.5ex>[r]^f \ar@<-0.5ex>[r]_g & B \ar@/^1.5pc/[l]^t \ar[r]^e & C\ar@/^/[l]^s , }
$$
such that $ef=eg,\, es=1,\, ft=1,\, gt=se$. A split coequaliser of $f$ and $g$ is the arrow $e$ of some split fork. 
\end{definition}

\begin{definition} \label{reflect}
A functor $U:\cD \to \cE$  is said to reflect coequalisers for a pair \mbox{$f,g:A \to B$} if for every diagram 
$\xymatrix@1{A \ar@<0.5ex>[r]^f \ar@<-0.5ex>[r]_g & B  \ar[r]^e & C , }$ with $ef=eg$ and $Ue$ a  coequaliser of $Uf$ and $Ug$, the morphism $e$ is a coequaliser of $f$ and $g$.

$U$ is said to reflect isomorphisms if $f$ is an isomorphism whenever $Uf$ is. Observe that if  $U$ preserves a coequaliser of $f,g$ in $\cD$, and 
$U$ reflects isomorphisms, then $U$ reflects all coequalisers of $f,g$.
\end{definition}

\begin{theorem}\label{beck} The following conditions are equivalent:
\begin{enumerate}

\item The adjunction 
$$
\xymatrix@1{\cD \ar@<1ex>[r]^U_{\top} & \cE \ar@<1ex>[l]^F},
$$
is monadic.

\item If $f,g: D \to D'$ is a parallel pair in $\cD$ for which $Uf,Ug$ has a split coequaliser, then $\cD$ has a coequaliser for $f,g$, and $U$ preserves and reflects coequalisers for these pairs.
\end{enumerate}
\begin{proof}
\cite{Mac} Ch. VI.7 Ex 6.
\end{proof}
\end{theorem}

Throughout this paper, the sign ``='' will be used to denote canonical isomorphism.

\section{Simplicial Deformation Complexes}\label{sdc}

\subsection{Definitions}
\begin{definition}\label{sdcdef} A simplicial deformation complex $E^{\bullet}$ consists of smooth homogeneous functors $E^n:\C_{\L} \to \Set$ for each $n \ge 0$, together with maps 
$$
\begin{matrix}
\pd^i:E^n \to E^{n+1} & 1\le i \le n\\
\sigma^i:E^{n}\to E^{n-1} &0 \le i <n,
\end{matrix}
$$
an associative product $*:E^m \by E^n \to E^{m+n}$, with identity $1: \bullet \to E^0$, where $\bullet$ is the constant functor $\bullet(A)=\bullet$ on $\C_{\L}$, such that:
\begin{enumerate}
\item $\pd^j\pd^i=\pd^i\pd^{j-1}\quad i<j$.
\item $\sigma^j\sigma^i=\sigma^i\sigma^{j+1} \quad i \le j$.
\item 
$
\sigma^j\pd^i=\left\{\begin{matrix}
			\pd^i\sigma^{j-1} & i<j \\
			\id		& i=j,\,i=j+1 \\
			\pd^{i-1}\sigma^j & i >j+1
			\end{matrix} \right. .
$
\item $\pd^i(e)*f=\pd^i(e*f)$.
\item $e*\pd^i(f)=\pd^{i+m}(e*f)$, for $e \in E^m$.
\item $\sigma^i(e)*f=\sigma^i(e*f)$.
\item $e*\sigma^i(f)=\sigma^{i+m}(e*f)$, for $e \in E^m$.
\end{enumerate}
\end{definition}

\begin{remark}
Observe that $E^{\bullet}$ is not quite a cosimplicial complex --- we have omitted the maps $\pd^0$ and $\pd^{n+1}$. In fact, if we define $\breve{E}_{\bullet-1}=E^{\bullet}$, with $\sigma_i=\pd^{i+1}$ and $\pd_{i}=\sigma^i$, then $\breve{E}_{\bullet}$ is an augmented simplicial complex.
\end{remark}

\begin{definition} Let $\CC^{\bullet}(E)$ be the tangent space of $E^{\bullet}$, i.e. $\CC^n(E)=E^n(k[\eps])$.
\end{definition}

\begin{definition}
Define the Maurer-Cartan functor $\mc_E$ by 
$$
\mc_E(A)=\{\omega \in E^1(A) \colon \omega*\omega=\pd^1(\omega)\}.
$$
\end{definition}

\begin{lemma}
$E^0$ is a group under multiplication.
\begin{proof}
We need only show existence of inverses. Given $g \in E^0(A)$, we have a morphism $\phi_g:E^0 \to E^0$ of functors on $\C_A$, given by $h \mapsto g*h$. Since $E^0$ is smooth, its obstruction theory is trivial. Now, on tangent spaces, $\phi_g (a)=\bar{g}*a$, but the map $A \to k[\eps]$ factors through $k$, so $\bar{g}=1$, and $\phi_g:\CC^0(E)\to \CC^0(E)$ is the identity. Thus the Standard Smoothness Criterion (Theorem \ref{SSC}) implies that $\phi_g$ is smooth, and in particular surjective, hence $\phi_g:E^0(A) \to E^0(A)$ is surjective, and there exists an $h$ with $g*h=1$. We construct a left inverse similarly, and then $h_l=h_l*g*h_r=h_r$. This also gives uniqueness of the inverse.
\end{proof}
\end{lemma}

Now, if $\omega \in \mc_E(A)$ and $g \in E^0(A)$, then $g*\omega*g^{-1}\in \mc_E(A)$. We may  therefore make the following definition:
\begin{definition}
$$
\ddef_E=\mc_E/E^0,
$$
the quotient being with respect to the adjoint action. The deformation groupoid
$$
\Def_E
$$
has objects $\mc_E$, and morphisms given by $E^0$.
\end{definition}

\begin{lemma}\label{faithful}
The action $E^0 \by E^n \to E^n$ is faithful (i.e. $g*e=e$ for some $e$ only if $g=1$).
\begin{proof}
From Theorem \ref{Man2}, it will suffice to show that $\nu:\CC(E)^0 \to \CC(E)^n$ is injective, since a deformation functor $F$ with tangent space $0$ has at most one element in $F(A)$.

Let $\omega_0$ be the unique element of $E^1(k)$; thus $\overbrace{\omega_0*\ldots*\omega_0}^n$ is the unique element of $E^n(k)$. Note that $\sigma^0(\omega_0)=1$.
Now, for $g \in \CC^0(E)$,
$$
(\sigma^0)^n(\nu(g))=(\sigma^0)^n(g*\omega_0^n)=(\sigma^0)^{n-1}(g*\sigma^0(\omega_0)*\omega_0^{n-1})=(\sigma^0)^{n-1}(g*\omega_0^{n-1})=g,
$$
by induction. Therefore $\nu$ is injective.
\end{proof}
\end{lemma}

Observe that, for $\omega \in \mc_E(A)$, 
$$
\sigma^0(\omega)*\omega=\sigma^0(\omega*\omega)=\sigma^0\pd^1(\omega)=\omega,
$$
so the previous Lemma implies that $\sigma^0(\omega)=1$. As an immediate consequence, we have:

\begin{lemma}
For $\omega \in \mc_E(A)$, if we let $\pd_{\omega}^0(e)=\omega*e$, and $\pd_{\omega}^{n+1}(e)=e*\omega$, for $e \in E^n(A)$, then $E^*(A)$ becomes a cosimplicial complex, $\pd_{\omega}^0$ and $\pd_{\omega}^{n+1}$ also satisfying  the following identities, for $e \in E^m, f \in E^n$:
\begin{eqnarray*}
\pd_{\omega}^0(e)*f&=&\pd_{\omega}^0(e*f),\\
e * \pd_{\omega}^{n+1}(f) &=& \pd_{\omega}^{m+n+1}(e*f),\\
\pd_{\omega}^{m+1}(e)*f&=& e* \pd_{\omega}^0(f).
\end{eqnarray*}
Note, however, that not all the identities of Definition \ref{sdcdef} hold for $\pd_{\omega}^0$ and $\pd_{\omega}^{n+1}$ --- in general $\pd_{\omega}^{m+1}(e)*f \ne \pd^{m+1}(e*f)$, and $e*\pd_{\omega}^0(f) \ne \pd^m(e*f)$.
\end{lemma}

In particular, we obtain, on $\C_k \subset \C_{\L}$, $\pd^0:=\pd^0_{\omega_0}$ and $\pd^{n+1}:=\pd^{n+1}_{\omega_0}$. Thus $\CC^{\bullet}(E)$ has the natural structure of a cosimplicial complex. 

\begin{definition}
Define the cohomology groups of $E$ to be
$$
\H^i(E):=\H^i(\CC^{\bullet}(E)),
$$
the cohomotopy groups of the cosimplicial complex $\CC^{\bullet}(E)$.
\end{definition}

\begin{lemma}
The tangent space of $\mc_E$ is $\z^1(\CC^{\bullet}(E))$, with the action of $E^0$ giving $\nu(g)=\pd^1(g)-\pd^0(g)$.
\begin{proof}
On tangent spaces all maps are linear. Thus the map $\CC^1(E)\by\CC^1(E)\to\CC^2(E)$ given by $(\alpha,\beta)\mapsto \alpha*\beta$ can be written $\alpha*\beta=f_1(\alpha)+f_2(\beta)$. Evaluating this at $\beta=\omega_0$ we obtain $f_1(\alpha)=\alpha*\omega_0=\pd^2(\alpha)$, and similarly $f_2(\beta)=\pd^0(\beta)$. Thus
$$
t_{\mc_E}=\{\omega \in \CC^1(E) \colon \pd^2(\omega)+\pd^0(\omega)=\pd^1(\omega)\}=\z^1(\CC^{\bullet}(E)).
$$

Now, $\nu(g)=g*\omega_0*g^{-1}$, and, similarly to above, 
$$
f*\omega_0*h=f*\omega_0*1+1*\omega_0*h=\pd^1(f)+\pd^0(h).
$$ 
Hence $\nu(g)=\pd^1(g)-\pd^0(g)$.
\end{proof}
\end{lemma}

\begin{lemma}\label{obstrsdc}
$\H^2(E)$ is a complete obstruction space for $\mc_E$.
\begin{proof}
Take a small extension
$$
e:0\to k \xra{\eps} A \to B \to 0.
$$
Note that $\eps $ provides an isomorphism $A\by_k k[\eps] \to A\by_B A$, $(\alpha,a\eps)\mapsto (\alpha, \alpha+a\eps)$, and thus, for any homogeneous functor $F$, $F(A)\by t_F\eps \to F(A\by_B A)$. Given $f \in F(A)$ and $v \in t_F$, we will write the image under this isomorphism as $(f,f+v\eps)$. Note that all the properties which this terminology suggests do hold.

Now, given $\omega \in \mc_E(B)$, use the smoothness of $E^1$ to take a lift $\tilde{\omega} \in E^1(A)$. Since $\omega*\omega=\pd^1(\omega)$, 
$$
(\tilde{\omega}*\tilde{\omega}, \pd^1(\tilde{\omega})) \in E^2(A)\by_{E^2(B)}E^2(A),
$$
 so we have $\tilde{\omega}*\tilde{\omega}=\pd^1(\tilde{\omega})+c\eps$, for $c \in \CC^2(E)$.

Now,
\begin{eqnarray*}
(\tilde{\omega}*\tilde{\omega})*\tilde{\omega}&=&\pd^1(\tilde{\omega})*\tilde{\omega}+\pd^3(c)\eps\\
\tilde{\omega}*(\tilde{\omega}*\tilde{\omega})&=&\tilde{\omega}*\pd^1(\tilde{\omega})+\pd^0(c)\eps,
\end{eqnarray*}
so
$$
\pd^1(\tilde{\omega}*\tilde{\omega})-\pd^2(\tilde{\omega}*\tilde{\omega})=(\pd^0(c)-\pd^3(c))\eps.
$$
Also,
\begin{eqnarray*}
\pd^1(\tilde{\omega}*\tilde{\omega})&=&\pd^1\pd^1(\tilde{\omega})+\pd^1(c)\eps\\
\pd^2(\tilde{\omega}*\tilde{\omega})&=&\pd^2\pd^1(\tilde{\omega})+\pd^2(c)\eps,
\end{eqnarray*}
so
$$
\pd^1(\tilde{\omega}*\tilde{\omega})-\pd^2(\tilde{\omega}*\tilde{\omega})=(\pd^1(c)-\pd^2(c))\eps,
$$
since $\pd^1\pd^1=\pd^2\pd^1$.

Thus $dc=(\pd^0-\pd^1+\pd^2-\pd^3)(c)=0$, so $c \in \z^2(\CC^{\bullet}(E))$.

If we take another lift $\tilde{\omega}_1=\tilde{\omega}+b\eps$, then
$$
\tilde{\omega}_1*\tilde{\omega}_1=\tilde{\omega}*\tilde{\omega}+\pd^0(b)\eps +\pd^2(b)\eps,
$$
since $\tilde{\omega}*b=\omega_0*b+\pd^0(b)$, {\&}c., and
$$
\pd^1(\tilde{\omega}_1)=\pd^1(\tilde{\omega})+\pd^1(b)\eps,
$$
so $c_1=c +db$.

We therefore define $v_e(\omega)=[c]\eps \in \H^2(E)$, and it follows from the above working that this is a complete obstruction theory.
\end{proof}
\end{lemma}

Theorems \ref{SSC} to \ref{keyhgs} now imply:

\begin{theorem}
$\ddef_E$ is a deformation functor, with tangent space $\H^1(E)$ and complete obstruction space $\H^2(E)$.  For $\omega, \omega' \in \mc_E(A)$, $\Iso(\omega,\omega')$ is homogeneous, with tangent space $\H^0(E)$ and complete obstruction space $\H^1(E)$. Moreover, if $\H^0(E)=0$, then $\ddef_E$ is homogeneous.
\begin{proof}
Theorem \ref{Man1} and Corollary \ref{keyhgs}.
\end{proof}
\end{theorem}

\begin{theorem}\label{sdcqis}
If $\phi :E \ra F$ is a morphism of SDCs, and 
$$
\H^i(\phi):\H^i(E) \ra \H^i(F)
$$ 
are the induced maps on cohomology, then:
\begin{enumerate}
\item If $\H^1(\phi)$ is bijective, and $\H^2(\phi)$ injective, then $\ddef_E \ra \ddef_F$ is \'etale.
\item If also $\H^0(\phi)$ is surjective, then  $\ddef_E \ra \ddef_F$ is an isomorphism.
\item Provided condition 1 holds, $\mathfrak{Def}_E \to \mathfrak{Def}_F$ is an equivalence of functors of  groupoids if and only if $\H^0(\phi)$ is an isomorphism.
\end{enumerate}
\begin{proof}
$ $
\begin{enumerate}
\item Theorems \ref{SSC} and \ref{Man1}.
\item Theorem \ref{Man3}.
\item It remains only to show that $\phi^0:E^0 \to F^0$ is injective. $E^0$ acts on $F^0$ by multiplication, so $\ker \phi^0 = \Iso(1,1)$, which, by Theorem \ref{Man2}, is a deformation functor with tangent space $\ker \H^0(\phi)=0$. Hence $\ker \phi^0=0.$ 
\end{enumerate}
\end{proof}
\end{theorem}

Call a morphism $\phi:E \to F$ a quasi-isomorphism if the $\H^i(\phi):\H^i(E) \ra \H^i(F)$ are all isomorphisms.

\begin{definition}\label{defmor} Given a morphism $\phi:E \to F$ of SDCs, and a point $x \in \mc_F(\L)$, define the groupoid
$$
\Def_{\phi,x}
$$
to be the fibre of the morphism
$$
\Def_E \to \Def_F
$$
over $x$.

Explicitly, $\Def_{\phi,x}(A)$ has objects 
$$
\{(\omega, h) \in \mc_E(A) \by F^0(A) \quad\colon\quad h\phi(\omega)h^{-1}=x\},
$$
and morphisms 
$$
E^0(A),\quad\text{ where }\quad  g(\omega,h)=(g\omega g^{-1},h\phi(g)^{-1}).
$$
\end{definition}

\begin{theorem}\label{cone}
Let $E^n_{\phi,x}(A)$ be the fibre of $E^n(A) \xra{\phi} F^n(A)$ over $x^n$. Then $E^{\bullet}_x$ is an SDC, and the canonical map
$$
\Def_{E_{\phi,x}} \to \Def_{\phi,x}
$$
is an equivalence of groupoids.
\begin{proof}
 $\ddef_{\phi,x}$ is a deformation functor, with tangent space $\H^1(E,F)$ and complete obstruction space $\H^2(E,F)$. For $\omega, \omega' \in \Ob \Def_{\phi,x}(A)$, $\Iso(\omega,\omega')$ is homogeneous, with tangent space $\H^0(E,F)$ and complete obstruction space $\H^1(E,F)$. Here, $\H^i(E,F)$ is the cohomology of the cone complex
$$
\xymatrix{
\vdots &\vdots\\
\CC^3(E) \ar[u]^d \ar[ur]^{-\phi} \ar@{}[r]|{\bigoplus}& \CC^2(F) \ar[u]^d\\
\CC^2(E) \ar[u]^d \ar[ur]^{\phi} \ar@{}[r]|{\bigoplus}& \CC^1(F) \ar[u]^d\\
\CC^1(E) \ar[u]^d \ar[ur]^{-\phi}\ar@{}[r]|{\bigoplus} & \CC^0(F) \ar[u]^d\\
\CC^0(E) \ar[u]^d \ar[ur]^{\phi}, 
}
$$
so we have the long exact sequence
$$
\xymatrix@R=0pt{
0 \ar[r]& \H^0(E,F)\ar[r]& \H^0(E)\ar[r]^{\H^0(\phi)}& \H^0(F)\ar[r]&\\
& \H^1(E,F)\ar[r]& \H^1(E)\ar[r]^{\H^1(\phi)}& \H^1(F)\ar[r]&\\
& \H^2(E,F)\ar[r]& \H^2(E)\ar[r]^{\H^2(\phi)}& \H^2(F)\ar[r]& \ldots.
}
$$

Now, the definition of $E^{\bullet}_x$ provides a short exact sequence:
$$
0 \to \CC^{\bullet}(E_{\phi,x}) \to \CC^{\bullet}(E) \xra{\phi} \CC^{\bullet}(F) \to 0,
$$
giving the long exact sequence
$$
\xymatrix@R=0pt{
0 \ar[r]& \H^0(E_{\phi,x})\ar[r]& \H^0(E)\ar[r]^{\H^0(\phi)}& \H^0(F)\ar[r]&\\
& \H^1(E_{\phi,x})\ar[r]& \H^1(E)\ar[r]^{\H^1(\phi)}& \H^1(F)\ar[r]&\\
& \H^2(E_{\phi,x})\ar[r]& \H^2(E)\ar[r]^{\H^2(\phi)}& \H^2(F)\ar[r]& \ldots.
}
$$
Thus the map $\Def_{E_{\phi,x}} \to \Def_{\phi,x}$ gives isomorphisms on all relevant tangent and obstruction spaces, so is an equivalence of groupoids.
\end{proof}
\end{theorem}

\subsection{Examples}

\subsubsection{Deformation of algebras}\label{alg}

Let $\L_n:=\L/\mu^{n+1}$.
Take a flat $\mu$-adic system of algebras $B_n/\L_n$, i.e  flat algebras $B_n/\L_n$ such that $B_{n+1}\ten_{\L_{n=1}}\L_n = B_n$, and let $B:=\Lim B_n$. The simplest (and most common) case to consider is when $B_n=\L_n$.  We wish to create an SDC describing flat deformations $R_A/(B\ten_{\L} A)$ of any flat algebra $R_0/B_0$:
$$
\begin{CD}
X_0=\Spec R_0 @>i>> \Spec R_A= X_A\\
 @V{\text{flat}}VV  @VV{\text{flat}}V  \\
  \Spec B_0 @>>> \Spf B \by_{\Spf\L} \Spec A.  
\end{CD}
$$
Explicitly, the deformation groupoid has objects flat $R_A/(B\ten A)$, such that \mbox{$R_A \ten_A k =R_0$,}  with morphisms given by infinitesimal isomorphisms of $R_A$, i.e. $(B\ten A)$-algebra isomorphisms $\theta$ such that \mbox{$\theta \equiv \id \mod \m_A$.}

There is, up to isomorphism, a unique flat $\mu$-adic $\L$-module $M$ such that $M\ten k= R_0$ (since flat modules over nilpotent rings are free).
There is an adjunction
$$
\xymatrix@=8ex{(B\ten A)\FAlg \ar@<1ex>[r]^-{U}_-{\top}	& A\FMod \ar@<1ex>[l]^-{F_{B \ten A}},}
$$
between the categories of flat $(B\ten A)$-algebras and  of  flat $A$-modules, $U$ being the forgetful functor, and 
$$
F_{B\ten A}=(B\ten A)\ten_A \Symm_{A}=B\ten_{\L} \Symm_{A}.
$$ 
Let $\bot=FU$ and $\top=UF$. We have the unit and co-unit:
$$
\eta_M: M = \rm{S}^1(M) \into  \Symm_{A}M \to B\ten_{\L}\Symm_{A}M = \top M,
$$
and
\begin{eqnarray*}
\vareps_R: \bot R &\to& R\\
b [r_1\ten\ldots\ten r_n]& \mapsto& b\cdot r_1\ldots  r_n.
\end{eqnarray*}
As in \cite{W}, we may form the canonical augmented cosimplicial complex
$$
\xymatrix@1{M \ar[r] &\top M \ar@<1ex>[r] \ar@<-1ex>[r]  &\top^2 M \ar@{.>}[l] \ar[r] \ar@/^/@<0.5ex>[r] \ar@/_/@<-0.5ex>[r] 
&\top^3M \ar@{.>}@<0.75ex>[l] \ar@{.>}@<-0.75ex>[l] \ar@/^1pc/[r] \ar@/_1pc/[r] \ar@{}[r]|{\cdot} \ar@{}@<1ex>[r]|{\cdot} \ar@{}@<-1ex>[r]|{\cdot} &\top^4 M & \ldots\ldots}
$$
with face operators 
$$
\top^i\eta\top^{n-i}:\top^n M \to \top^{n+1}M
$$ 
and degeneracy operators 
$$
\top^iU\vareps F\top^{n-i}:\top^{n+1}M \to \top^{n}M.
$$

\begin{definition}
Let 
$$
E^n=\Hom_{\L}(\top^n M,M)_{U(\vareps_{R_0}\circ\ldots\circ\vareps_{\bot^{n-1} R_0})},
$$
where, for flat $\mu$-adic $\L$-modules $M$ and $N$, and a $k$-linear map $\alpha:M\ten k \to N\ten k$, we define 
$$
\Hom_{\L}(M,N)_{\alpha}(A)
$$ 
to be the fibre of 
$$
\Hom_{A}(M\ten A,N\ten A) \to \Hom_{k}(M\ten k, N\ten k)
$$ 
over $\alpha$.

We make $E^*$ into an SDC by giving it the product $g*h=g\circ\top^n h$, and $\breve{E}_*$ is the dual complex to the canonical augmented complex $M \to \top^*M$.
\end{definition}

Since flat modules over Artinian rings do not deform (being free), to give a deformation $R/(B\ten_{\L}A)$ of $R_0/B_0$ is the same as to give the module $M\ten_{\L}A$ the structure of a $(B\ten_{\L}A)$-algebra, compatible with the $B_0$-algebra structure of $R_0$, and the $A$-module structure of $M\ten_{\L}A$.

 A $(B\ten_{\L}A)$-algebra structure on  module $M_A:=M\ten_{\L} A$ is the same as a map
$$
f \in \Hom_{A\FMod}(\top M_A, M_A); \quad b[ m_1\ten\ldots\ten m_n] \mapsto b\cdot m_1\ldots  m_n,
$$
satisfying
\begin{eqnarray*}
f\circ\top f &=& f \circ U\vareps_{FM_A},\\
 f\circ \eta_{M_A} &=& \id_{M_A}.
\end{eqnarray*}

Hence elements of $\mc_E(A)$ correspond to  \mbox{$(B\ten A)$}-algebra structures on $M\ten A$, compatible with $R_0$ (recall that, for $f \in \mc_E$, we automatically have $\sigma^{0}(f)=1$). Elements $f_1,f_2$ will give isomorphic ring structures precisely when we have $g \in E^0(A)$ such that
$$
\begin{CD}
	\top M \ten A 	@>{f_1}>> M \ten A	\\
@V{\top g}VV 			@VV{g}V\\
	\top M	\ten A	@>{f_2}>> M\ten A	
\end{CD}
$$
commutes.

Therefore $\Def_E \simeq \Def_{R_0/B}$.

\begin{remark}
Observe that, if $\omega \in \mc_E(A)$ gives rise to a ring $R$, then we have an isomorphism of cosimplicial complexes
$$
(E^n(A),\pd_{\omega})\cong \Hom_{B\ten A\Alg}(\bot_{*}R,R)_{\vareps_R \circ \vareps_{\bot R}\circ \ldots \circ \vareps_{\bot^n R}},
$$
where $\bot_n R =\bot^{n+1}R$. The natural product $(f,g) \mapsto f \circ \bot^{n+1}g$, for $f:\bot_n R \to R$, corresponds to the product $(\alpha,\beta) \mapsto \alpha * \omega * \beta$ in $E^*(A)$.

In particular, there is an isomorphism of cosimplicial complexes
$$
\CC^{\bullet}(E) \cong \Der_{B_0}(\bot_{\bullet}R_0,R_0) = \Hom_{R_0}(\mathbf{L}_{\bullet}^{R_0/B_0}, R_0),
$$
where $\mathbf{L}_{\bullet}^{R_0/B_0}$ is one presentation of the cotangent complex. Illusie (\cite{Ill1} and \cite{Ill2}) and Andr\'e (\cite{An}) define the cotangent complex using the resolution arising from the adjunction between algebras and sets, whereas we have used the adjunction between flat algebras and flat modules. However, Quillen (\cite{Q}) proves that all simplicial resolutions give the same object in the derived category.  

Therefore $\H^*(E)$ is Andr\'e-Quillen cohomology.
\end{remark}

If $R_0/B_0$ is smooth, and lifts to some ($\mu$-adic) smooth algebra $R/B$, then we have a morphism of SDCs
$$
\Hom_{B\FAlg}(R,R)_{\id} \to E^n,
$$
where the complex on the left has the same object at every level, the obvious multiplication map, and all face and degeneracy operators given by the identity. Since $R/B$ is smooth, this is a quasi-isomorphism, so the first complex also describes deformations of $R_0$.

\subsubsection{Deformation of sheaves --- \v Cech resolution}\label{zarshf}

Take a  sheaf of $k$-vector spaces $\sM_0$ on a topological space $X$. The deformation functor will associate to $A$, flat sheaves $\sM_A$ of $A$-modules such that $\sM_A\ten k=\sM_0$, modulo infinitesimal isomorphisms.

Take an open cover $\mathfrak{U}$ of $X$ 
such that $U_{i_1}\cap U_{i_2}\cap \ldots \cap U_{i_n}$ is contractible for all $i_1, \ldots i_n$. Examples for which this is possible include $X$ a manifold or a geometric simplicial complex.
Set 
$$
X'=\coprod_i U_i.
$$
We have inclusion maps $u_i: U_i \to X$, so may combine them to give  $u:X' \to X$. 

Now, since flat modules do not deform, there are, up to isomorphism,  unique $\mu$-adic $\L$-modules $\sN_i$ lifting $\Gamma(U_i,\sM_0)$. Combining these, we obtain a sheaf $\sN$ on $X'$ with  \mbox{$\sN\ten k = u^{-1}\sM_0$,} and there is an adjunction
$$
\xymatrix@1{\Shf(X')\ar@<1ex>[r]^{u_*}_{\top}	&\Shf(X) \ar@<1ex>[l]^{u^{-1}}},
$$
with unit
\begin{eqnarray*}
\alpha:\sF &\to& u_*u^{-1}\sF,\\
\Gamma(U,\sF) &\to& \Gamma(U,u_*u^{-1}\sF)=\prod_i \Gamma(U\cap U_i,\sF)\\
f &\mapsto& \prod_i \rho_{U,U\cap U_i}(f),
\end{eqnarray*}
and co-unit
\begin{eqnarray*}
\gamma:  u^{-1}u_*\sG &\to& \sG,\\
\prod_j (\sG_j)|_{U_{ji}}= \Gamma(U_i,u^{-1}u_*\sG) &\to& \sG_i \\
g &\mapsto& g_i \in \sG_{ii}=\sG_i.
\end{eqnarray*}

\begin{definition}
We define the SDC $E^*$ by:
$$
E^n=\Hom_{\L}(\sN, (u^{-1}u_*)^n(\sN))_{\alpha^n},
$$
where $\alpha^n=u^{-1}(\alpha_{(u_*u^{-1})^{n-1}\sM_0}\circ \ldots \circ \alpha_{\sM_0})$, with the SDC structure arising exactly as in the previous section (modulo contravariance).
\end{definition}

Now, to give a sheaf $\sM$ of $A$-modules on $X$, such that $u^{-1}\sM = \sN\ten_{\L}A$, is to give transition functions
$$
\theta_{ji}: \sN_i\ten A \to (\sN_j\ten A)|_{U_{ij}},
$$
satisfying
$$
\theta_{kj}\circ\theta_{ji}=\theta_{ki};\quad \theta_{ii}=\id.
$$

This is precisely the same as a map
$$
\Theta: \sN \ten A \to u^{-1}u_*( \sN \ten A),
$$ 
satisfying
$$
(u^{-1}u_*\Theta) \circ \Theta = u^{-1}\alpha_{u_*\sN_A}\circ \Theta;\quad \gamma_{\sN}\circ \Theta = \id_{\sN}.
$$

If we include the additional condition that these transition functions must be compatible with those of $\sM_0$, mod $ \m_A$, then we have precisely an element of $\mc_E(A)$. It is easy to see that equivalences of transition functions correspond to elements of $E^0(A)$, acting via the adjoint action.

Note that we may regard $E^n$ as 
$$
\Gamma(X', \hom_{\L}(\sN, (u^{-1}u_*)^n(\sN))_{\alpha^n}),
$$ 
and that setting 
$$
\sE^n=u_*\hom_{\L}(\sN, (u^{-1}u_*)^n(\sN))_{\alpha^n},
$$
 with tangent space $\CC^n(\sE)$, we have that the chain complex $\CC^{\bullet}(\sE)$  is a flabby resolution of the sheaf $\hom_k(\sM_0,\sM_0)$ (it is, in fact,a version of  the \v Cech resolution), so that
$$
\H^i(E)=\H^i(\CC^{\bullet}(E))=\bH^i(X,\CC^{\bullet}(\sE))=\H^i(X,\hom_k(\sM_0,\sM_0)),
$$
as expected. Equivalently,
$$
\H^i(\CC^{\bullet}(E))=\Ext^i_k(\sM_0,\sM_0)=\H^i(X,\hom_k(\sM_0,\sM_0)),
$$
since we are working with sheaves of vector spaces over $k$, so $\hom_k$ is exact.

\section{A Procedure for Computing SDCs}

In this section we make formal the approach which has been used so far to compute SDCs. The idea is that we throw away properties of the object which we wish to deform, until we obtain something whose deformations are trivial.

Throughout this section, we will encounter functors $\cD:\C_{\L} \to \Cat$. We will \emph{not} require that these functors satisfy (H0) (the condition that $F(k)=\bullet$).

\begin{definition} Given a functor $\cD:\C_{\L} \to \Cat$, and an object $D \in \Ob\cD(k)$, define
$\Def_{\cD,D}:\C_{\L} \to \Grpd$ by setting $\Def_{\cD,D}(A)$ to be the fibre of $\cD(A) \to \cD(k)$ over $(D,\id)$.
\end{definition}
All the deformation problems we encounter are of this form.

\begin{definition}
We say a functor $\cB:\C_{\L} \to \Cat$ has uniformly trivial deformation theory if the functor $\Mor\cB$ is smooth and homogeneous, and the functor $\Cmpts\cB$ is  constant,  i.e for $A \onto A'$ in $\C_{\L}$, $\cB(A) \to \cB(A')$ is full and essentially surjective. A typical example of such a functor is that which sends $A$ to the category of flat $A$-modules.

Observe that, if $\cB$ is uniformly trivial, then given $B \in \Ob \cB(k)$, we may lift it to $\tilde{B} \in \Ob \cB(\L)$, and we have an equivalence of functors of groupoids
$$
(\tilde{B}, \Mor_{\cB}(\tilde{B},\tilde{B})) \xra{\sim} \Def_{\cB,B}.
$$
\end{definition}

\subsection{The SDC associated to a monadic adjunction}

Assume we have a monadic adjunction
$$
\xymatrix{\cD \ar@<1ex>[r]^{U}_{\top} &\cB \ar@<1ex>[l]^F},
$$
with unit $\eta:\id \to UF=\top$ and co-unit $\vareps:FU \to \id$, such that $\cB$ has uniformly trivial deformation theory. Then, given $D \in \cD(k)$, let $B$ be any lift of $UD$ to $\cB(\L)$, and define the SDC $E^*$ by
$$
E^n=\Hom_{\cB}(\top^n B,B)_{U(\vareps_{D}\circ\ldots\circ\vareps_{\bot^{n-1} D})},
$$
the fibre of
$$
\Hom_{\cB}(\top^n B,B) \to \Hom_{\cB(k)}(\top^n B(k), B(k))
$$
over $U(\vareps_{D}\circ\ldots\circ\vareps_{\bot^{n-1} D})$.

We make $E^*$ into an SDC by giving it the product $g*h=g\circ\top^n h$, and $\breve{E}_*$ is the dual complex to the canonical augmented complex $B \to \top^*B$. Explicitly, for $g \in E^n$,
\begin{eqnarray*}
\pd^i(g) &=& g \circ \top^{i-1}U\vareps_{F\top^{n-i}B}\\
\sigma^i(g) &=& g \circ \top^{i}\eta_{\top^{n-i-1}B}.
\end{eqnarray*}

\begin{theorem}
With the notation as above, we have an equivalence of functors of groupoids
$$
\Def_{\cD,D} \xra{\sim} \Def_E.
$$
\begin{proof}
Since the $F \dashv U$ is monadic, we have an equivalence $K: \cD \to \cB^{\top}$. Now observe that $\mc_E$ is precisely the fibre of
$$
\Ob \cB^{\top} \to \Ob \cB^{\top}(k)
$$
over $K(D)$ (note the comment after Lemma \ref{faithful} gives $\sigma^0(\omega)=1$), with morphisms in the fibre
$$
\cB^{\top} \to \cB^{\top}(k)
$$
over $(K(D), \id)$ corresponding to $E^0$ acting on $\mc_E$ via the adjoint action.
Hence
$$
\Def_E= \Def_{\cB^{\top},K(D)}.
$$
But $K: \cD \to \cB^{\top}$ is an equivalence, giving an equivalence
$$
K:\Def_{\cD,D} \xra{\sim}  \Def_{\cB^{\top},K(D)}=  \Def_E.
$$
\end{proof}
\end{theorem}

Observe that this result is precisely the same as the description above of deformations of algebras, and that the dual result (for comonadic adjunctions) describes deformations of sheaves.

\subsubsection{Deformation of sheaves --- Godement resolution}\label{etshf}
Take a sheaf of $k$-vector spaces $\sM_0$ on a site $X$ with enough points. In examples, we will be considering the Zariski or \'etale sites of a scheme, or the site associated to a topological space. The deformation functor will associate to $A$,  sheaves $\sM_A$ of flat $A$-modules such that $\sM_A\ten k=\sM_0$, modulo infinitesimal isomorphisms.

Let $X'$ be the set of points of $X$.
Since $X$ has enough points, the inverse image functor $\Shf(X) \to \prod_{x \in X'} \Shf(x)$ reflects isomorphisms. Explicitly, this says that a morphism $\theta: \sF \to \sG$ of sheaves on $X$ is an isomorphism whenever the morphisms $\theta_x:\sF_x \to \sG_x$  are for all $x \in X'$. In the reasoning which follows, it will suffice to replace $X'$ by any subset with this property. 

If we are working on the Zariski site, we may take $X'=\coprod_{x \in  X} x$. On the \'etale site, we may take
$$
X'=\coprod_{x \in  X} \bar{x},
$$
where for each $x \in X$, a geometric point $\bar{x} \to X$ has been chosen. For Jacobson schemes, we may consider only closed points $x \in X$. We define the category $\Shf(X')$ by $\Shf(X'):=\prod_{x \in X'} \Shf(x)$. 

There is an  adjunction
$$
\xymatrix@1{\Shf(X')\ar@<1ex>[r]^{u_*}_{\top}	&\Shf(X) \ar@<1ex>[l]^{u^{*}}},
$$
the maps $u_x: x \to X$ combining to form  maps \mbox{$u^*:\Shf(X) \to \Shf(X')$,} and \mbox{$u_*: \Shf(X')\to \Shf(X)$} given by $u_*\sF=\prod_{x \in X'} u_{x*}\sF_x$. Observe that the category of flat $A$-modules on $X'$ has uniformly trivial deformation theory.

It follows from Theorem \ref{beck} that this adjunction is comonadic, and from Theorem \ref{main}, the  SDC governing this problem is
$$
E^n=\Hom_{\L}(\sN, (u^*u_*)^n(\sN))_{\alpha^n},
$$
where  $\sN$ is a flat $\mu$-adic $\L$-module   on $X'$ with $\sN\ten k = u^*\sM_0$, and  \mbox{$\alpha^n=u^*(\alpha_{(u_*u^*)^{n-1}\sM_0}\circ \ldots \circ \alpha_{\sM_0}).$}

To see this more clearly in the case of a topological space, observe, as in \cite{Go}, Remark II.4.3.2, that an element $n \in \Gamma(U,u^{-1}u_*\sN)$ 
can be represented by a function of the form
$$
n(x_0,x_1) \in \sN_{x_1},
$$
defined on a set of the form
$$
    x_0 \in U;\,         x_1 \in U(x_0),
$$
where $U(x_0)$ is some open \'etale neighbourhood of $x_0$ in $X$. Two functions $n_1, n_2$ define the same section if they agree on some set of the above form. 

 Thus an element $\Theta$ of $\mc_E(A)$ can be regarded as local transition functions; given $n \in \sN_x$, we obtain, for some \'etale neighbourhood $U$ of $x$, $\Theta(n)_{xy} \in \sN_y$ for all $y \in U$, which can be regarded as lifting the germ and then projecting down to each stalk. The Maurer-Cartan equation is transitivity, and, as before, guarantees that $\Theta(n)_{xx}=n$. 

Again, we get
$$
\H^i(E)=\H^i(\CC^{\bullet}(E))=\bH^i(X,\CC^{\bullet}(\sE))=\H^i(X,\hom_k(\sM_0,\sM_0))=\Ext^i_k(\sM_0,\sM_0),
$$
as expected.

\subsubsection{Deformation of co-algebras}\label{coalg}

Given a flat (co-associative) co-algebra (with co-unit) $C_0/k$,  we wish to create an SDC describing flat deformations $C_A/A$ such that $C_A \ten_A k = C_0$, modulo infinitesimal isomorphisms.

There is, up to isomorphism, a unique flat $\mu$-adic $\L$-module $M$ such that \mbox{$M\ten_{\L} k= C_0$}. 
There is an adjunction
$$
\xymatrix@1{ A\mathrm{-FCoAlg} \ar@<1ex>[r]^V & \ar@<1ex>[l]^G_{\bot} A \FMod},
$$
between the category of  flat co-algebras over $A$, and the category of  flat modules over $A$,
where $V$ is the forgetful functor and the free functor $G$ exists by the Special Adjoint Functor Theorem (since $A\mathrm{-FCoAlg}$ has all colimits, and $V$ preserves these). Note that in this case the free functor is hard to write down explicitly, but this is unnecessary for our purposes. See \cite{Sweedler} for such a description.

By the Theorem \ref{beck}, this adjunction is comonadic,   so that deformations of $C_0$ are given by the SDC
$$
E^n=\Hom_{\L}(M,\bot^n M)_{U(\alpha_{\top^{n-1} C_0}\circ\ldots\circ\alpha_{C_0})}.
$$

\subsection{The general approach}\label{gensdc}

In general, we will not be able to pass from the category $\cD$ to a category $\cB$ with uniformly trivial deformation theory via a single  monadic or comonadic adjunction. However, we should be able to pass from $\cD$ to some $\cB$ via a chain of monadic and comonadic adjunctions.

Not only should we have  monadic and comonadic adjunctions, but, informally, the forgetful functors should commute with one another. More precisely, assume that we have a diagram
$$
\xymatrix@=8ex{
\cD \ar@<1ex>[r]^{U}_{\top} \ar@<-1ex>[d]_{V} 
&\ar@<1ex>[l]^{F} \cE  \ar@<-1ex>[d]_{V} 
\\
\ar@<-1ex>[u]_{G}^{\dashv}	\cA \ar@<1ex>[r]^{U}_{\top} 
&\ar@<1ex>[l]^{F} \ar@<-1ex>[u]_{G}^{\dashv} \cB,  
}
$$
where $\cB$ has uniformly trivial deformation theory,
with $F\dashv U$ monadic and  $G\vdash V$ comonadic. Let 
\begin{align*}
\toph&=UF&	\both&=FU\\
\botv&=VG&	\topv&=GV,
\end{align*}
with 
$$
\eta:1 \to \toph,\quad \gamma:\botv \to 1,\quad \vareps:\both \to 1 \text{ and }\alpha:1 \to \topv.
$$
The commutativity condition is that the following identities  hold:
\begin{equation}
GU=UG \quad \text{ or } \quad  FV=VF  \label{eqnone}
\end{equation}
\begin{equation}
UV=VU
\end{equation}
\begin{equation}
V\vareps=\vareps V \quad	\text{ or }\quad U\alpha=\alpha U\\
\end{equation}
\begin{equation}
V\eta=\eta V	\quad	\text{ or }\quad	U\gamma=\gamma U. \label{eqnfive}
\end{equation}
Observe that the adjoint properties ensure that identities on the same line are equivalent.

\begin{lemma}\label{delta}
For $E \in \cE,  A \in \cA$, consider the diagram of isomorphisms given by the adjunctions:
$$
\xymatrix@=2ex{
&\Hom_{\cD}(FE,GA) \ar@{<->}[dr] \ar@{<->}[dl]& \\
 \Hom_{\cA}(VFE,A)\ar@{=}[d] &	&\Hom_{\cE}(E,UGA)\ar@{=}[d]\\
\Hom_{\cA}(FVE,A)\ar@{<->}[dr]	&	&\Hom_{\cE}(E,GUA)\ar@{<->}[dl]\\
&\Hom_{\cB}(VE,UA). &}
$$
The identities (\ref{eqnone})--(\ref{eqnfive}) ensure that all squares in this diagram commute. In particular, for any $B \in \cB$, the map 
$$
\eta_{B}\circ\gamma_{B}: VGB \to UFB
$$
corresponds to a map
$$
\rho_B=\vareps_{GFB}\circ FG\eta_B = GF\gamma_B \circ \alpha_{FGB}:FGB \to GFB,
$$
with
\begin{eqnarray*}
\gamma_{\toph B}\circ UV(\rho_B)\circ \eta_{\botv B}&=& \eta_{B}\circ\gamma_{B},\\
U(\rho_B)\circ \eta_{GB}&=& G(\eta_B),\\
\gamma_{FB}\circ V(\rho_B) &=& F(\gamma_B).
\end{eqnarray*}
\end{lemma}

By the naturality of $\vareps:\both \to \id$, we have $\vareps \circ (\both \vareps) =\vareps \circ (\vareps \both)$, so we obtain a canonical map $\vareps^n: \both^n \to \id$, given by any such composition of $\vareps$'s. [For instance 
$$
\vareps^n=\vareps\circ (\both\vareps)\circ\ldots \circ (\both^{n-1}\vareps)
$$
is one such composition.]
 We have similar maps for each of the units and co-units, giving:
\begin{align*}
\vareps^n:\both^n &\to \id,	&\eta^n:\id &\to \toph^n \\
\alpha^n:\id &\to \topv^n,	&\gamma^n:\botv^n &\to \id 
\end{align*}

Let 
$$
\delta = UV(\rho):\toph\botv \to \botv \toph.
$$
Since $\delta$ is natural, 
$$
(\delta \botv\toph)\circ(\toph\botv \delta)=(\botv\toph \delta)\circ (\delta \toph\botv).
$$
Therefore any composition of $\delta$'s gives us the same canonical map
$$
\delta^{m,n}:\toph^m\botv^n \to \botv^n\toph^m.
$$

\begin{theorem}\label{main}
Suppose we have a diagram as above. Then, for $D \in \Ob\cD(k)$, let $B$ be any lift of $UV D \in \Ob\cB(k)$ to $\Ob\cB(\L)$ (by the hypothesis on $\cB$, such a lift must exist and be unique up to isomorphism). Set
$$
E^n=\Hom_{\cB}(\toph^n B, \botv^n B)_{UV(\alpha^n_{D}\circ \vareps^n_D)}.
$$
We give $E^*$ the multiplication 
$$
g*h= \botv^n(g)\circ\delta^{m,n}\circ \toph^m(h),
$$ 
and $\breve{E}_*$ gets its simplicial structure from the diagonal of the canonical double augmented complex 
$$
\Hom_{\cB}(B \to  B^*, B_* \to B).
$$
Explicitly, for $g \in E^n$,
\begin{eqnarray*}
\pd^i(g) &=&    \botv^{i-1}V\alpha_{G\botv^{n-i}B}\circ	g \circ \toph^{i-1}U\vareps_{F\toph^{n-i}B}\\
\sigma^i(g) &=&	 \botv^{i}\gamma_{\botv^{n-i-1}B}\circ	g \circ \toph^{i}\eta_{\toph^{n-i-1}B}.
\end{eqnarray*}

Then 
$$
\Def_{\cD,D} \simeq \Def_{E}.
$$
\begin{proof}
First observe that this is, indeed, an SDC. We need to check that the product is associative. For $f \in \SDC^l, g \in \SDC^m, h \in \SDC^n$, we have
\begin{eqnarray*}
f*(g*h)&=& \botv^{m+n}f \circ \delta^{l,m+n}\circ\toph^l(\botv^n g \circ \delta^{m,n}\circ\toph^m h)\\
	&=& \botv^{m+n}f \circ (\botv^n\delta^{l,m}\circ \delta^{l,n}\botv^m)\circ \toph^l\botv^n g  \circ\toph^l\delta^{m,n}\circ\toph^{l+m} h\\
	&=& \botv^{m+n}f \circ \botv^n\delta^{l,m}\circ\botv^n\toph^l g \circ \delta^{l,n}\toph^m \circ\toph^l\delta^{m,n}\circ\toph^{l+m} h \\
	&=& \botv^n(\botv^m f \circ \delta^{l,m}\circ \toph^l g) \circ \delta^{l+m,n}\circ \toph^{l+m} h\\
	&=& (f*g)*h.
\end{eqnarray*}

We now have a lemma:
\begin{lemma}
To give a map $f:\toph B \to \botv B$ satisfying the Maurer-Cartan equation
$$
\botv f \circ \delta \circ \toph f = V \alpha_{GB} \circ f \circ U\vareps_{FB}
$$
is the same as giving maps $\theta: \toph B \to B$ and $\phi: B \to \botv B$ satisfying the Maurer-Cartan equations:
\begin{eqnarray*}
\theta \circ \toph \theta &=& \theta \circ U\vareps_{FB}\\
\botv \phi \circ\phi  &=& V \alpha_{GB} \circ\phi,
\end{eqnarray*}
and the compatibility condition
$$
\botv\theta \circ \delta_B \circ\toph \phi = \phi \circ \theta.
$$
\begin{proof}\emph{(Of Lemma.)}

Given such an $f$, we have
\begin{eqnarray*}
f &=&\gamma_{\botv B} \circ V \alpha_{GB} \circ f \circ U\vareps_{FB} \circ \eta_{\toph B}\\
&=& \gamma_{\botv B} \circ \botv f \circ\delta_B\circ \toph f \circ \eta_{\toph B}\\
&=& f \circ \gamma_{\toph B} \circ\delta_B\circ \eta_{\botv B} \circ f \\
&=& f \circ \eta_B \circ \gamma_B \circ f.
\end{eqnarray*}
Let $\phi=f \circ \eta_B$, and $\theta=\gamma_B \circ f$.
Now,
\begin{eqnarray*}
\gamma_B \circ \gamma_{\botv B} \circ \botv f \circ\delta_B\circ \toph f &=& \gamma_B \circ \gamma_{\botv B}\circ V \alpha_{GB} \circ f \circ U\vareps_{FB}\\
\gamma_B \circ f \circ \gamma_{\toph B}\circ\delta_B\circ \toph f &=& \gamma_B \circ f \circ U\vareps_{FB}\\
\gamma_B \circ f \circ U(\gamma_{F B}\circ V(\rho_B))\circ \toph f &=& \gamma_B \circ f \circ U\vareps_{FB}\\
\gamma_B \circ f \circ \toph\gamma_B\circ \toph f &=& \gamma_B \circ f \circ U\vareps_{FB}\\
\theta\circ \toph \theta &=& \theta \circ U\vareps_{FB},
\end{eqnarray*}
using the Maurer-Cartan equations and Lemma \ref{delta}. Similarly we obtain the Maurer-Cartan equation for $\phi$, and finally
\begin{eqnarray*}
\botv\theta \circ\delta_B\circ \toph \phi&=&\botv\gamma_B \circ \botv f\circ\delta_B \circ \toph f \circ \toph\eta_B\\
&=&\botv\gamma_B \circ V \alpha_{GB} \circ f \circ U\vareps_{FB} \circ \toph\eta_B\\
&=& f\\
&=& \phi \circ \theta.
\end{eqnarray*}

Conversely, given such $\theta$ and $\phi$, set $f=\phi \circ \theta$. We obtain:
\begin{eqnarray*}
\botv f \circ\delta_B\circ \toph f &=& \botv\phi \circ \botv \theta\circ\delta_B \circ \toph\phi  \circ \toph\theta\\
&=& \botv \phi \circ\phi \circ \theta \circ \toph \theta\\
&=& V \alpha_{GB} \circ\phi \circ \theta \circ U\vareps_{FB}\\
&=& V \alpha_{GB} \circ f \circ U\vareps_{FB}.
\end{eqnarray*}
\end{proof}
\end{lemma}

Recall, from the comment after Lemma \ref{faithful}, that every element $f$ of $\mc_E(A)$ must satisfy $\sigma^0(f)=1$, which in this case is $\gamma_B \circ f \circ \eta_B = \id$. This means that the pair $(\theta, \phi)$ above satisfy $\theta \circ \eta_B = \id$, and $ \gamma_B \circ \phi = \id$.

We obtain an adjunction
$$
\xymatrix@=8ex{\cE^{\toph} \ar@<1ex>[r]^-{V} & \cB^{\toph} \ar@<1ex>[l]^-{G}_-{\bot}},
$$
\begin{eqnarray*}
\left(\vcenter{\xymatrix{
\toph E \ar[d]_f \\ E }}\right) &\mapsto& \left(\vcenter{\xymatrix{\toph V E \ar[d]_{Vf} \\ VE }}\right)\\
\left(\vcenter{\xymatrix{
\toph G B \ar[d]_{G(g)\circ U(\rho_B)}\\ GB}}\right) &\mapsfrom& \left(\vcenter{\xymatrix{\toph B \ar[d]_{g} \\ B }}\right),
\end{eqnarray*}
and such a pair $(\theta, \phi)$ is precisely the same as an element of $(\cB^{\toph})^{\botv}$. We also have a correspondence between the action of $E^0$ on $\mc_{E}$ and morphisms in $(\cB^{\toph})^{\botv}$,    so that $\Def_{E}$ corresponds to the fibre of
$$
(\cB^{\toph})^{\botv} \to (\cB^{\toph})^{\botv}(k)
$$
over the image of $(D,\id)$.

Now, the equivalence
\begin{eqnarray*}
K: \cA &\to& \cB^{\toph}\\
A &\mapsto& \left(\vcenter{\xymatrix{ \toph UA  \ar[d]_{U \vareps_A} \\
				 UA}} \right)
\end{eqnarray*}
arising from the monadic adjunction $F \dashv U$ satisfies
\begin{eqnarray*}
\botv K A &=& VG \left(\vcenter{\xymatrix{ \toph U A \ar[d]_{ U \vareps_{A}} \\  UA}}\right)\\
	&=&\left(\vcenter{\xymatrix{\toph\botv U A \ar[d]_{\botv U \vareps_{A}\circ VU(\rho_{UA})}  \\  \botv UA}}\right).
\end{eqnarray*}
Moreover, 
\begin{eqnarray*}
\botv U \vareps_{A}\circ VU(\rho_{UA}) &=& UV(G\vareps_A \circ (\vareps_{GFUA}\circ FG\eta_{UA}))\\
&=&UV((\vareps_{GA}\circ FUG\vareps_{A})\circ FG\eta_{UA})\\
&=&UV(\vareps_{GA} \circ FG(U\vareps_A \circ \eta_{UA}))\\
&=& UV(\vareps_{GA}),
\end{eqnarray*}
so
$$
\botv K A = \left(\vcenter{\xymatrix{\toph\botv U A \ar[d]_{ UV(\vareps_{GA})} \\  \botv UA}}\right) = K \botv A.
$$
This induces an equivalence of categories
$$
 \cA^{\botv} \xra{K_{\mathrm{h}}}(\cB^{\toph})^{\botv}. 
$$
But $G \vdash V$ is comonadic, so 
$$
\cD \xra{K_{\mathrm{v}}}\cA^{\botv}
$$
is an equivalence.
Thus we have an equivalence of groupoids
$$
\Def_{\cD,D} \xra{\sim}\Def_{\cA^{\botv}, K_{\mathrm{v}} D} \xra{\sim} \Def_{(\cB^{\toph})^{\botv},K_{\mathrm{h}} K_{\mathrm{v}} D} =    \Def_E.
$$
\end{proof}
\end{theorem}

Given a category $\cD$, we now need a procedure for  finding a category $\cB$ with uniformly trivial deformation theory. We do this by thinking of the functors $U$ and $V$ as being forgetful functors, and looking for structure in $\cD$ to discard.

\begin{definition}
Given a category $\cB:\C_{\L} \to \Cat$, define a set of structures $\Sigma$ over $\cB$ to consist of the following data:
\begin{enumerate}
\item A finite set $\Sigma= \Sigma^+\sqcup \Sigma^-$,
\item For each  subset $S \subset \Sigma$, a category $\cB^S$, with $\cB^{\emptyset}=\cB$,
\item \begin{enumerate} 
	\item for each $s \in \Sigma^+$ and each $S \subset \Sigma$ not containing $s$,  a monadic adjunction
	$$
	\xymatrix@1{
	\cB^{S\cup \{s\}} \ar@<1ex>[r]^-{U_s}_-{\top} &\ar@<1ex>[l]^-{F_s} \cB^S}.
	$$
	\item for each $s \in \Sigma^-$ and each $S\subset \Sigma$ not containing $s$,  a comonadic adjunction
	$$
	\xymatrix@1{
	\cB^{S\cup \{s\}} \ar@<1ex>[r]^-{V_s} &\ar@<1ex>[l]^-{G_s}_-{\bot} \cB^S},
	$$
	\end{enumerate}
\end{enumerate}
satisfying the  commutativity conditions given below.

For $s \in \Sigma^+$ we write
$$
\top_s=U_s F_s \quad \text{ and } \quad	\bot_s=F_s U_s,
$$
with 
$$
\eta_s:1 \to \top_s \quad \text{ and }\quad \vareps_s:\bot_s \to 1.
$$
For $s \in \Sigma^-$ we write
$$
\bot_s=V_s G_s	\quad \text{ and } \quad \top_s=G_s V_s,
$$
with 
$$
\gamma_s:\bot_s \to 1,\quad \text{ and } \quad \alpha_s:1 \to \top_s.
$$

The commutativity conditions are that the following  hold:
\begin{enumerate}
\item For each distinct pair $s_1,s_2 \in \Sigma^+$, and each $S\subset \Sigma$ not containing $s_1,s_2$, the morphisms in the diagram
$$
\xymatrix@=8ex{
\cB^{S\cup \{s_1,s_2\}} \ar@<1ex>[r]^-{U_1}_-{\top} \ar@<-1ex>[d]_{U_2}^{\vdash} 
&\ar@<1ex>[l]^-{F_1} \cB^{S\cup \{s_2\}}  \ar@<-1ex>[d]_{U_2}^{\vdash} 
\\
\ar@<-1ex>[u]_{F_2}	\cB^{S\cup \{s_1\}} \ar@<1ex>[r]^-{U_1}_-{\top} 
&\ar@<1ex>[l]^-{F_1} \ar@<-1ex>[u]_{F_2} \cB^S  
}
$$
satisfy
\begin{equation*}
F_1F_2=F_2F_1 \quad \text{ or } \quad  U_1U_2=U_2U_1  
\end{equation*}
\begin{align*}
U_1\vareps_2&=\vareps_2 U_1&	U_2\vareps_1&=\vareps_1 U_2\\
U_1\eta_2&=\eta_2 U_1&			U_2\eta_1&=\eta_1 U_2. 
\end{align*}

\item For each distinct pair $s_1,s_2 \in \Sigma^-$, and each $S \subset \Sigma$ not containing $ s_1,s_2$, the morphisms in the diagram
$$
\xymatrix@=8ex{
\cB^{S\cup \{s_1,s_2\}} \ar@<1ex>[r]^-{V_1} \ar@<-1ex>[d]_{V_2} 
&\ar@<1ex>[l]^-{G_1}_-{\bot} \cB^{S\cup \{s_2\}}  \ar@<-1ex>[d]_{V_2} 
\\
\ar@<-1ex>[u]_{G_2}^{\dashv}	\cB^{S\cup \{s_1\}} \ar@<1ex>[r]^-{V_1} 
&\ar@<1ex>[l]^-{G_1}_-{\bot} \ar@<-1ex>[u]_{G_2}^{\dashv} \cB^S  
}
$$
satisfy
\begin{equation*}
G_1G_2=G_2G_1 \quad \text{ or } \quad  V_1V_2=V_2V_1  
\end{equation*}
\begin{align*}
V_1\alpha_2&=\alpha_2 V_1&	V_2\alpha_1&=\alpha_1 V_2\\
V_1\gamma_2&=\gamma_2 V_1&			V_2\gamma_1&=\gamma_1 V_2. 
\end{align*}

\item For each $s_1 \in \Sigma^+$ and $s_2 \in \Sigma^-$, and  each $S \subset \Sigma$ not containing $ s_1,s_2$, the morphisms in the diagram
$$
\xymatrix@=8ex{
\cB^{S\cup \{s_1,s_2\}} \ar@<1ex>[r]^{U}_{\top} \ar@<-1ex>[d]_{V} 
&\ar@<1ex>[l]^{F} \cB^{S\cup \{s_2\}}  \ar@<-1ex>[d]_{V} 
\\
\ar@<-1ex>[u]_{G}^{\dashv}	\cB^{S\cup \{s_1\}} \ar@<1ex>[r]^{U}_{\top} 
&\ar@<1ex>[l]^{F} \ar@<-1ex>[u]_{G}^{\dashv} \cB^S,  
}
$$
satisfy
\begin{equation*}
GU=UG \quad \text{ or } \quad  FV=VF  
\end{equation*}
\begin{equation*}
UV=VU
\end{equation*}
\begin{equation*}
V\vareps=\vareps V \quad	\text{ or }\quad U\alpha=\alpha U\\
\end{equation*}
\begin{equation*}
V\eta=\eta V	\quad	\text{ or }\quad	U\gamma=\gamma U. 
\end{equation*}
\end{enumerate}
\end{definition}

\begin{lemma} Given a set of structures $\Sigma$ over a category $\cB$, we have a diagram
$$
\xymatrix@=8ex{
\cB^{\Sigma} \ar@<1ex>[r]^{U}_{\top} \ar@<-1ex>[d]_{V} 
&\ar@<1ex>[l]^{F} \cB^{\Sigma^-}  \ar@<-1ex>[d]_{V} 
\\
\ar@<-1ex>[u]_{G}^{\dashv}	\cB^{\Sigma^+} \ar@<1ex>[r]^{U}_{\top} 
&\ar@<1ex>[l]^{F} \ar@<-1ex>[u]_{G}^{\dashv} \cB,  
}
$$
satisfying the  equations (\ref{eqnone})--(\ref{eqnfive}) on p. \pageref{eqnone},  with $F\dashv U$ monadic and  $G\vdash V$ comonadic.

\begin{proof} We define $F$ (resp. $U$) to be the composition of the $F_s$ (resp $U_s$) for all $s \in \Sigma^+$, noting that the order of composition does not matter, since these functors commute with one another. We define $G$ and $V$ analogously. It is immediate that $F \dashv U$ and $G \vdash V$ are adjoint pairs, and that the commutativity conditions (\ref{eqnone})--(\ref{eqnfive}) are satisfied. 

It thus remains only to show that $F \dashv U$ is monadic, and $G \vdash V$ comonadic. This can be done by using a similar approach to the proof of Theorem \ref{main}. The statement is that $\cB^{\Sigma^+} \simeq \cB^{\top_+}$, where $\top_+=FU$, which is proved   by induction on the cardinality of $\Sigma^+$.
\end{proof}
\end{lemma}

\begin{remarks}
\begin{enumerate}
\item This lemma shows that  every set of structures can be replaced by a set of at most two elements, so it might seem that introducing the notion of a set of structures was not helpful. However, it is frequently easier to find the individual adjunctions than the composites, as will be seen in later examples (notably the deformation of a group scheme in Section \ref{gpscheme}).

\item When considering those deformation problems for which DGLAs have been successfully constructed, it seems that  there is a  set of structures  which is either wholly monadic or wholly comonadic.
\end{enumerate}
\end{remarks}

\subsubsection{Deformation of Hopf algebras}

Given a flat (associative, commutative, co-associative) Hopf algebra (with unit and co-unit) $R_0/k$,  we wish to create an SDC describing flat deformations $R_A/A$ of $R_0$ such that $R_A \ten_A k = R_0$, modulo infinitesimal isomorphisms.

The structures are:
$$
\Sigma^+=\{\text{Algebra}\},\quad \Sigma^-=\{ \text{Co-Algebra}\},
$$
over the category  of flat $A$-modules.

This gives the following commutative diagram of monadic and comonadic adjunctions:
$$
\xymatrix@C=12ex@R=8ex{
A\mathrm{-FHopfAlg} \ar@<1ex>[r]_-{\top} \ar@<-1ex>[d]
& \ar@<1ex>[l]^-{\Symm_{A}}	A\mathrm{-FCoAlg}\ar@<-1ex>[d]
\\
\ar@<-1ex>[u]_{G}^{\dashv} A\FAlg \ar@<1ex>[r]_-{\top}
&\ar@<-1ex>[u]_{G }^{\dashv} \ar@<1ex>[l]^-{\Symm_{A}} A\FMod,}
$$
where 
$A\mathrm{-FHopfAlg}$ is the category of  flat Hopf algebras over $A$, $A\mathrm{-FCoAlg}$ is the category of  flat co-associative co-algebras with co-unit and co-inverse over $A$, and
$G$ is the free co-algebra functor of Section \ref{coalg}.  Since we cannot describe $G$  explicitly, we use the alternative form of the axioms involving $\Symm_{A}$ instead, when verifying the conditions of Theorem \ref{main}.

We thus obtain the SDC
$$
E^n=\Hom_{\L}((\Symm_{\L})^n M,  G^n M)_{\alpha^n \circ \vareps^n},
$$
where $M$ is the flat $\mu$-adic $\L$-module (unique up to isomorphism) lifting the   $k$-vector space $R_0$.
with $\alpha^n, \vareps^n$  the canonical maps 
$$
\vareps^n:(\Symm_k)^n R_0 \to R_0,\quad
\alpha^n: R_0  \to G^n R_0
$$
associated to $R_0$.   

\section{Deformation of Schemes}\label{scheme}

Given a flat $\mu$-adic system of schemes $S_n/\Spec \L_n$, and a flat scheme $X_0/S_0$, our deformation functor consists of  flat  schemes $X_A/S_A$ such that \mbox{$X_A\by_{\Spec A}\Spec k =X_0$,} modulo infinitesimal isomorphisms (isomorphisms which pull back to the identity on $X_0$),  where $\fS=\varinjlim S_n$ and $S_A=\fS\by_{\Spf\L}\Spec A$:
$$
\begin{CD}
X_0				@>>>	X_A		\\
@V{ f \text{ flat}}VV	 		@VV{\text{flat}}V\\
S_0 				@>>>	  S_A.
\end{CD}
$$

Since the topological space $|X_0|$ underlying $X_0$ does not deform, this is just a question of deforming the sheaf $\O_{X_0}$ of algebras. Observe that, since $X_0$ is flat over $S_0$, any scheme $X_A/S_A$ deforming $X$, and flat over $A$, will necessarily be flat over $S_A$. It therefore suffices to deform $\O_{X_0}$ as an $\O_{S_A}$-algebra, flat over $A$.

In the notation of Section \ref{gensdc}, we have the set of structures 
$$
\Sigma^+=\{\text{Algebra},\,\, \O_{S_A}\text{-Module}\},\quad \Sigma^-=\{X_0\text{-Sheaf}\},
$$
over the category  
of sheaves of flat $A$-modules on $X_0'$, where $X_0'$ is defined as in Section \ref{etshf}.
To these structures correspond the monadic and comonadic adjunctions
$$
\{\Symm_{A}{}\dashv U,\,\,\O_{\fS}\ten_{\L}{} \dashv U\},\quad \{ u_* \vdash u^{-1}\},
$$
where $U$ denotes the relevant forgetful functor.

This yields the following diagram of $\Cat$-valued functors: 
$$
\xymatrix@C=12ex@R=8ex{
\O_{S_A}\backslash(A\FAlg(X_0)) \ar@<1ex>[r]_{\top} \ar@<-1ex>[d]_{u^{-1}}
&\ar@<1ex>[l]^{\O_{\fS}\ten_{\L}\Symm_{A}} A\FMod(X_0)  \ar@<-1ex>[d]_{u^{-1}}
\\
\ar@<-1ex>[u]_{u_*}^{\dashv}	\O_{S_A}\backslash(A\FAlg(X_0')) \ar@<1ex>[r]_{\top} 
&\ar@<1ex>[l]^{\O_{\fS}\ten_{\L}\Symm_{A}} \ar@<-1ex>[u]_{u_*}^{\dashv} A\FMod(X_0'). 
}
$$
where we write $\O_{\fS}$ for the sheaf $f^{-1}\O_{\fS}$ (resp. $u^{-1}f^{-1}\O_{\fS}$) on $X_0$ (resp. $X'_0$), and $\O_{S_A}\backslash(A\FAlg)$ consists  of those  $\O_{Y_A}$-algebras which are flat over $A$. The only non-trivial commutativity condition is the observation that pull-backs commute with tensor operations.

Hence, by Theorem \ref{main}, deformations are described by the SDC
$$
E^n=\Hom_{\L}((\O_{\fS}\ten_{\L}\Symm_{\L})^n \sM, (u^{-1}u_*)^n \sM)_{u^{-1}(\alpha^n\circ \vareps^n)},
$$
where $\sM$ is a lift of the sheaf $u^{-1}\O_{X_0}$ of  vector spaces on $X_0'$ to a sheaf of flat $\mu$-adic $\L$-modules, and $\alpha^n,\vareps^n$ are the canonical maps $\alpha^n_{\O_{X_0}},\vareps^n_{\O_{X_0}}$ in  Theorem \ref{main}, given by the adjunctions.  

\begin{remark}
Observe that the description above allows us to construct an SDC governing deformations of a sheaf of algebras on any site with enough points. This is more general than the problem in \cite{Hinich} for which a DGLA was constructed.
\end{remark}

Define a sheaf of SDCs on $X$ by 
$$
\sE^{n}=u_*\hom_{\L}((\O_{\fS}\ten_{\L}\Symm_{\L})^n \sM, (u^{-1}u_*)^n \sM)_{u^{-1}(\alpha^n\circ \vareps^m)}.
$$
Combining the observations in  Sections \ref{alg}, \ref{etshf}, we see that this has tangent space 
\begin{eqnarray*}
\CC^{n}(\sE)&=& u_*\der_{\O_{S_0}}((\O_{S_0}\ten_{k}\Symm_{k})^{n+1} u^{-1}\O_{X_0}, (u^{-1}u_*)^n u^{-1}\O_{X_0})\\
&\cong& \der_{\O_{S_0}}((\O_{S_0}\ten_{k}\Symm_{k})^{n+1} \O_{X_0}, (u_* u^{-1})^{n+1}\O_{X_0})\\
&\cong& \hom_{\O_{X_0}}(\mathbf{L}_n^{X_0/S_0},\sC^n(\O_{X_0})),
\end{eqnarray*}
where $\mathbf{L}_{\bullet}^{X_0/S_0}$ is the standard form of the cotangent complex, as described in \cite{Ill1}, and $\sC^n$ denotes the Godement resolution. 

Since $C^{n}(\sE)$ has the structure of a diagonal complex of a bicosimplicial complex, it follows from the Eilenberg-Zilber Theorem that it is quasi-isomorphic to the total complex of the double complex
$$
\hom_{\O_{X_0}}(\mathbf{L}_{\bullet}^{X_0/S_0},\sC^{\bullet}(\O_{X_0}))
$$
Therefore the cohomology of our SDC is
$$
\H^i(\CC^{\bullet}(E))= \H^i(\Gamma(X_0, \CC^{\bullet}(\sE)))\cong  \EExt^i_{\O_{X_0}}(\mathbf{L}_{\bullet}^{X_0/S_0}, \O_{X_0}),
$$
 which is Andr\'e-Quillen hypercohomology, the second isomorphism following because $\mathbf{L}_{\bullet}^{X_0/S_0}$ is locally projective as an $\O_{X_0}$-module, and $\sC^{\bullet}(\O_{X_0})$ is flabby. 

\subsection{Separated Noetherian schemes}\label{sepn}
If $X_0$ is separated and Noetherian, then we may replace Godement resolutions by \v Cech resolutions. Take an open affine cover $(X_{\alpha})_{\alpha \in I}$ of $X_0$, and set $\check{X}:=\coprod_{\alpha \in I}X_{\alpha}$. We then have a diagram
$$
\xymatrix{\check{X}' \ar[r]^{\check{u}} \ar[dr]^w \ar[d]_{v'} & \check{X}\ar[d]^v \\
X_0' \ar[r]^u & X_0.}
$$
Since $v^{-1}\O_{X_0}=\O_{\check{X}}$ is a quasi-coherent sheaf on $\check{X}$, $\Ext^i(\O_{\check{X}},\O_{\check{X}})=0$ for all $i>0$, so deformations of the sheaf $\O_{\check{X}}$ of $k$-vector spaces on $\check{X}$ are unobstructed and we may lift $v^{-1}\O_{X_0}$ to some sheaf $\sN$ of flat $\mu$-adic $\L$-modules on $\check{X}$, unique up to non-unique isomorphism. There must also be (non-canonical) isomorphisms $\check{u}^{-1}\sN \cong {v'}^{-1}\sM$. If we define $\check{E}^{\bullet}$ by 
$$
\check{E}^n=\Hom_{\L}((\O_{\fS}\ten_{\L}\Symm_{\L})^n \sN, (v^{-1}v_*)^n \sN)_{v^{-1}(\alpha^n\circ \vareps^n)},
$$
and ${F}^{\bullet}$ similarly, replacing $v$ by $w$, then we have morphisms of SDCs
$$
E^{\bullet} \xra{\gamma_{v'}^{\bullet}} F^{\bullet} \xleftarrow{\gamma_{\check{u}}^{\bullet}}\check{E}^{\bullet}.
$$
To see that these are quasi-isomorphisms, let $\sT^{\bullet}:=\hom_{\O_{X_0}}(\mathbf{L}_{\bullet}^{X_0/S_0},\O_{X_0})$, which is a complex of quasi-coherent sheaves on $X_0$, and observe that the induced maps on cohomology are
$$
\bH^i(X_0, \sT^{\bullet}) \to \bH^i(X_0, \sT^{\bullet}) \leftarrow \check{\bH}^i(X_0, \sT^{\bullet}).
$$
The final map is an isomorphism since $X_0$ is separated and Noetherian, so \v Cech cohomology agrees with sheaf cohomology. 

\subsection{Smooth schemes}\label{smooth}
If $X_0/S_0$ is smooth, then the cohomology groups  are  just $\H^i(X_0,\sT_{X_0/S_0})$, the cohomology of the tangent sheaf. If $X_0$ is also separated and Noetherian (for instance if $S_0$ is so), then we may lift $\check{X}$ to some smooth $\mu$-adic formal  scheme $\check{\fX}$ over $\fS$.

Consider the diagram
$$
\xymatrix@1{\O_{S_A}\backslash(A\FAlg(\check{X}))\ar@<1ex>[r]^{v_*}_{\top}	&\O_{S_A}\backslash(A\FAlg(X_0)) \ar@<1ex>[l]^{v^{*}}};
$$
although the former category does not have uniformly trivial deformation theory, all the morphisms in it which we encounter do, which gives us an SDC and 
 canonical maps
$$
\Hom_{\O_{\fS}\Alg}(\O_{\check{\fX}}, (v^{-1}v_*)^n\O_{\check{\fX}})_{v^{-1}(\alpha^n)} \xra{(\vareps^n)^*} \check{E}^n,
$$
which give  a quasi-isomorphism of SDCs.

Explicitly, we may write our SDC as 
$$
E^n = \prod_{\alpha_0,\ldots, \alpha_n \in I} \Gamma(X_{\alpha_0,\ldots, \alpha_n}, \hom_{\O_{\fS}\Alg}(\sO_{\fX_{\alpha_n}},\sO_{\fX_{\alpha_0}})_{\id}),
$$
where $X_{\alpha_0,\ldots, \alpha_n}=\bigcap_{i=1}^n X_{\alpha_i}$, and $\check{\fX}=\coprod \fX_{\alpha}$. The product is given by
$$
(\phi*\psi)_{\alpha_0,\ldots, \alpha_{m+n}}= \phi_{\alpha_0,\ldots, \alpha_m}\circ \psi_{\alpha_m,\ldots, \alpha_{m+n}},
$$
and operations are 
\begin{eqnarray*}
(\pd^i\phi)_{\alpha_0,\ldots,\alpha_{n+1}}&=& \phi_{\alpha_0,\ldots, \alpha_{i-1},\alpha_{i+1},\ldots, \alpha_{n+1}}\\
(\sigma^i\phi)_{\alpha_0,\ldots,\alpha_{n-1}}&=& \phi_{\alpha_0,\ldots, \alpha_{i},\alpha_{i},\ldots, \alpha_{n-1}}
\end{eqnarray*}

\subsection{Deformation of a subscheme}\label{subscheme}

The standard application of the previous section would be to set $S_n=\Spec \L_n$. However, exactly the same proof applies to a far wider class of problems. Rather than requiring both $X_A/S_A$ and $\fS/\Spf\L$ to be flat, we only used the fact that $X_A/\Spec A$ was flat (flatness of $X_A/S_A$ then being automatic, since $X_0/S_0$ flat). This gives us the following immediate generalisation:

Given
$$
\xymatrix{X_0 \ar[dr]_{\text{flat}} \ar[r] &Y_0 \ar[d]^{\text{flat}}\\
					&S_0\ar[d] \\
&\Spec k,}
$$
for flat $\mu$-adic families
$$
\Y \xra{\text{flat}} \fS \xra{\text{flat}} \Spf \L,
$$
let $Y_A:=\Y\ten_{\L}A$, $S_A:=\fS\ten_{\L} A$, and  consider deformations
$$
\xymatrix{
X_0 \ar[dr] \ar[r] \ar[rrd] &Y_0 \ar[d] \ar[rrd]\\
&S_0 \ar[rrd] &X_A \ar[dr]|{\text{flat}} \ar[r] & Y_A \ar[d]^{\text{flat}}\\
&&&S_A,
}
$$
such that $X_A\ten_{A}k=X_0$.

Thus we need to consider  $\O_{Y_A}$-algebras deforming  $\O_{X_0}$, flat over $\O_{S_A}$. These are the same as $\O_{Y_A}$-algebras deforming  $\O_{X_0}$, flat over $A$. 

We have the set of structures
$$
\Sigma^+=\{\text{Algebra},\,\, \O_{Y_A}\text{-Module}\},\quad \Sigma^-=\{X_0\text{-Sheaf}\},
$$
over the category $A\FMod(X_0')$, giving  the diagram:
$$
\xymatrix@C=12ex@R=8ex{
\O_{Y_A}\backslash(A\FAlg(X_0)) \ar@<1ex>[r]_-{\top} \ar@<-1ex>[d]_{u^{-1}} 
&\ar@<1ex>[l]^{\O_{\fY}\ten_{\L}\Symm_{A}} A\FMod(X_0)  \ar@<-1ex>[d]_{u^{-1}}
\\
\ar@<-1ex>[u]_{u_*}^{\dashv}	\O_{Y_A}\backslash(A\FAlg(X_0')) \ar@<1ex>[r]_-{\top} 
&\ar@<1ex>[l]^{\O_{\fY}\ten_{\L}\Symm_{A}} \ar@<-1ex>[u]_{u_*}^{\dashv} A\FMod(X_0') 
}
$$
 and we obtain the same SDC as before, with $\fS$ replaced by $\Y$. Hence deformations are described by the SDC
$$
E^n=\Hom_{\L}((\O_{\Y}\ten_{\L}\Symm_{\L})^n \sM, (u^{-1}u_*)^n \sM)_{u^{-1}(\alpha^n\circ \vareps^n)},
$$
for a flat $\mu$-adic $\L$-module $\sM$ lifting $u^{-1}\O_{X_0}$.

In particular this applies to deformation of a subscheme $X_0 \into Y_0$, since a flat deformation of a subscheme as a $Y_A$-scheme will be a subscheme.  
If $X_0 \into Y_0$ is a regular embedding (e.g. if $X_0$ and $Y_0$ are both smooth over $S_0$), then the cohomology groups will be
$$
\H^i(\CC^{\bullet}(E))=\H^{i+1}(X_0, \sN_{X_0/Y_0}),
$$
where $\sN_{X_0/Y_0}$ is the normal sheaf.

If $X_0$ and $Y_0$ are both smooth over $S_0$, we can simplify the SDC of Section \ref{sepn} still further. Consider the diagram
$$
\xymatrix@C=12ex@R=8ex{
\O_{Y_A}\backslash(A\FAlg(X_0)) \ar@<1ex>[r]_{\top} \ar@<-1ex>[d]_{v^{-1}} 
&\ar@<1ex>[l]^{\O_{\fY}\ten_{\O_{\fS}}} \O_{S_A}\backslash(A\FAlg(X_0))  \ar@<-1ex>[d]_{v^{-1}}
\\
\ar@<-1ex>[u]_{v_*}^{\dashv}	\O_{Y_A}\backslash(A\FAlg(\check{X})) \ar@<1ex>[r]_{\top} 
&\ar@<1ex>[l]^{\O_{\fY}\ten_{\O_{\fS}}} \ar@<-1ex>[u]_{v_*}^{\dashv} \O_{S_A}\backslash(A\FAlg(\check{X})). 
}
$$
Although the last category does not have uniformly trivial deformation theory, all the morphisms in it which we encounter do.
If we let $(f^{\sharp})^{\ten n}$ denote the canonical ring homomorphism  $\O_{\fY}^{\ten_{\O_{\fS}} n}\ten_{\O_{\fS}} \O_{\hat{\fX}} \to \O_{\hat{\fX}}$, then there is an SDC and a quasi-isomorphism
$$
\Hom_{\O_{\fS}\Alg}(\O_{\fY}^{\ten_{\O_{\fS}} n}\ten_{\O_{\fS}} \O_{\hat{\fX}}, (v^{-1}v_*)^n\O_{\hat{\fX}})_{v^{-1}\alpha^n\circ (f^{\sharp})^{\ten n}} \to \check{E}^n.
$$  
The map is defined by composing the maps 
$
\Symm_{\L}(\O_{\fY}^{\ten n} \ten \O_{\hat{\fX}}) \to \O_{\fY}^{\ten n}\ten \O_{\hat{\fX}},
$
arising from the natural ring structure on the tensor product.

Explicitly, we may write our SDC as 
$$
E^n = \prod_{\alpha_0,\ldots, \alpha_n \in I} \Gamma(X_{\alpha_0,\ldots, \alpha_n}, \hom_{\O_{\fS}\Alg}(\O_{\fY}^{\ten_{\O_{\fS}} n}\ten_{\O_{\fS}} \O_{\fX_{\alpha_n}},\sO_{\fX_{\alpha_0}})_{(f^{\sharp})^{\ten n}}),
$$
The product is given by
\begin{eqnarray*}
&(\phi*\psi)_{\alpha_0,\ldots, \alpha_{m+n}}(r_1\ten \ldots\ten r_{m+n}\ten a)=&\\ &\phi_{\alpha_0,\ldots, \alpha_m}(r_1\ten \ldots\ten r_{m}\ten \psi_{\alpha_m,\ldots, \alpha_{m+n}}(r_{m+1}\ten \ldots\ten r_{m+n}\ten a))&,
\end{eqnarray*}
and operations are  
\begin{eqnarray*}
&&\quad (\pd^i\phi)_{\alpha_0,\ldots,\alpha_{n+1}}(r_1\ten \ldots\ten r_{n+1}\ten a)=\\ &&\phi_{\alpha_0,\ldots, \alpha_{i-1},\alpha_{i+1},\ldots, \alpha_{n+1}}(r_1\ten \ldots\ten r_{i-1}\ten (r_i r_{i+1})\ten r_{i+2}\ldots\ten r_{n+1}\ten a)\\
&&\quad (\sigma^i\phi)_{\alpha_0,\ldots,\alpha_{n-1}}(r_1\ten \ldots\ten r_{n-1}\ten a)=\\ &&\phi_{\alpha_0,\ldots, \alpha_{i},\alpha_{i},\ldots, \alpha_{n-1}}(r_1\ten \ldots\ten r_i\ten 1\ten r_{i+1}\ldots\ten r_{n-1}\ten a).
\end{eqnarray*}

To see that the map is a quasi-isomorphism, observe that the Eilenberg-Zilber Theorem allows us to regard the cohomology of $E$ as hypercohomology of the complex $\sT_{X_0/S_0} \to f^*\sT_{Y_0/S_0}[-1]$, which is quasi-isomorphic to the tangent complex of the morphism $X_0 \to Y_0$. Note that this SDC works for any morphism of smooth schemes, not only for embeddings. 

\subsection{Constrained deformation of schemes}

Given a morphism of flat $\mu$-adic schemes $\fZ \xra{h} \fY$ over $\fS$, and a diagram 
$$
Z_0 \xra{g} X_0 \xra{f} Y_0
$$ 
over $S_0$, with $fg=h$,  our deformation functor consists of  diagrams 
$$Z_A \xra{g_A} X_A \xra{f_A}Y_A,$$ 
over $S_A$, with $X_A \to S_A$ flat and $f_A \circ g_A=h$, which pull back along $S_0 \to S_A$ to our original diagram. Here, $Z_A$ is the scheme $\fZ\by_{\Spf \L} \Spec A$, and $Z_0=Z_k$. An example of such a deformation problem would be deformations of a subscheme, constrained to pass through a fixed set of points. 

Since the topological space $|X_0|$ underlying $X_0$ does not deform,  it suffices to deform $\O_{X_0}$ as an $\O_{Y_A}$-algebra with an $\O_{Z_A}$-augmentation, flat over $A$.

In the notation of Section \ref{gensdc}, we have the set of structures 
$$
\Sigma^+=\{\text{algebra},\,\, \O_{Y_A}\text{-module}\},\quad \Sigma^-=\{X_0\text{-sheaf},\O_{Z_A}\text{-augmented} \},
$$
over the category  
of sheaves of flat $A$-modules on $X_0'$, where $X_0'$ is defined as in Section \ref{etshf}.

To these structures correspond the monadic and comonadic adjunctions
$$
\{\Symm_{A}{}\dashv U,\,\,\O_{\fY}\ten_{\L}{} \dashv U\},\quad \{ u_* \vdash u^{-1},\O_{\fZ}\by{}\vdash V \},
$$
where $U,V$ denote the relevant forgetful functors.

This yields the following diagram of $\Cat$-valued functors: 
$$
\xymatrix@C=12ex@R=8ex{
\O_{Y_A}\backslash(A\FAlg(X_0))/\O_{Z_A} \ar@<1ex>[r]_-{\top} \ar@<-1ex>[d]_{u^{-1}} 
&\ar@<1ex>[l]^-{\O_{\fY}\ten_{\L}\Symm_{A}} A\FMod(X_0)/ \O_{Z_A}  \ar@<-1ex>[d]_{u^{-1}}
\\
\ar@<-1ex>[u]_{\O_{Z_A}\by u_*}^{\dashv}	\O_{Y_A}\backslash(A\FAlg(X_0')) \ar@<1ex>[r]_-{\top} 
&\ar@<1ex>[l]^-{\O_{\fY}\ten_{\L}\Symm_{A}} \ar@<-1ex>[u]_{\O_{Z_A}\by u_*}^{\dashv} A\FMod(X_0'). 
}
$$
where we write $\O_{\fY}$ for the sheaf $f^{-1}\O_{\fY}$ (resp. $u^{-1}f^{-1}\O_{\fY}$), and $\O_{\fZ}$ for the sheaf $g_*\O_{\fZ}$ (resp. $u^{-1}g_*\O_{\fZ}$)  on $X_0$ (resp. $X'_0$). The category $\O_{Y_A}\backslash(A\FAlg)/ \O_{Z_A}$ consists  of those  $\O_{Y_A}$-algebras which are flat over $A$, and equipped with an augmentation to $\O_{Z_A}$, compatible with $h^{\sharp}$. The only non-trivial commutativity condition is the observation that pull-backs commute with tensor operations.

By Theorem \ref{main}, this deformation problem is governed by the SDC
$$
E^n=\Hom_{\L}((\O_{\fY}\ten_{\L}\Symm_{\L})^n \sM, (\O_{\fZ}\by u^{-1}u_*)^n \sM)_{u^{-1}(\alpha^n\circ \vareps^n)},
$$
where $\sM$ is a lift of the sheaf $u^{-1}\O_{X_0}$ of  vector spaces on $X_0'$ to a sheaf of flat $\mu$-adic $\L$-modules, and $\alpha^n,\vareps^n$ are the canonical maps $\alpha^n_{\O_{X_0}},\vareps^n_{\O_{X_0}}$ in  Theorem \ref{main}, given by the adjunctions. The cohomology is 
$$
\EExt^i_{\O_{X_0}}(\mathbf{L}_{\bullet}^{X_0/Y_0}, \O_{X_0}\xra{g^{\sharp}}g_{*}\O_{Z_0}[-1] ),
$$

If $X_0$ and $Y_0$ are both smooth over $S_0$, we may consider the diagram
$$
\xymatrix@C=12ex@R=8ex{
\O_{Y_A}\backslash (A\FAlg(X_0))/ \O_{Z_A} \ar@<1ex>[r]_-{\top} \ar@<-1ex>[d]_{v^{-1}} 
&\ar@<1ex>[l]^-{\O_{\fY}\ten_{\O_{\fS}}} \O_{S_A}\backslash (A\FAlg(X_0))/ \O_{Z_A}  \ar@<-1ex>[d]_{v^{-1}}
\\
\ar@<-1ex>[u]_{ \O_{Z_A}\by v_*}^{\dashv}	\O_{Y_A}\backslash(A\FAlg(\check{X}))/ \O_{Z_A}  \ar@<1ex>[r]_-{\top} 
&\ar@<1ex>[l]^-{\O_{\fY}\ten_{\O_{S_A}}} \ar@<-1ex>[u]_{ \O_{Z_A}\by v_*}^{\dashv} \O_{S_A}\backslash (A\FAlg(\check{X})). 
}
$$
All the morphisms which we encounter in the last category  have uniformly trivial deformation theory, allowing us to replace this SDC by
$$
E^n = \prod_{\alpha_0,\ldots, \alpha_n \in I} \Gamma(X_{\alpha_0,\ldots, \alpha_n}, \hom_{\O_{\fS}\Alg}(\O_{\fY}^{\ten_{\O_{\fS}} n}\ten\O_{\fS} \O_{\fX_{\alpha_n}},\sO_{\fX_{\alpha_0}}\by \O_{\fZ}^n)_{(g^{\sharp})^n\circ(f^{\sharp})^{\ten n}}),
$$
since  cohomology in this case is hypercohomology of the complex 
$$
\sT_{X_0/S_0} \to (f^*\sT_{Y_0/S_0} \oplus \sT_{X_0/S_0}\ten_{\O_{X_0}}g_{*}\O_{Z_0})[-1]\to f^*\sT_{Y_0/S_0}\ten_{\O_{X_0}}g_{*}\O_{Z_0}[-2].
$$

\subsection{Deformation of a pair}
Given a flat $\mu$-adic formal scheme $\fS$, and  flat schemes 
$$
\xymatrix@=2ex{
X_0 \ar[rr]^{f_0} \ar[ddr] &&Y_0 \ar[ddl] 
\\
\\			&S_0,}
$$
let $S_A:=\fS\by_{\Spf\L}\Spec A$. Our deformation functor consists of (infinitesimal) isomorphism classes of:
$$
\xymatrix@=2ex{
X_0 \ar[rr]^{f_0} \ar[ddr] \ar[rrrrdd] 
& 
&Y_0  \ar[ddl] \ar[rrrrdd] 
\\
\\
&S_0 \ar[rrrrdd]
&&
&X_A \ar[rr]^{f_A} \ar[ddr]|{\text{flat}}
&
&Y_A  \ar[ddl]^{\text{flat}}
\\
\\
&&&&
& S_A,
}
$$
 such that $X_A\by_{\Spec A}\Spec k =X_0$,  $Y_A\by_{\Spec A}\Spec k =Y_0$, and $f_A\ten_A k =f_0$.

In this case, we obtain the SDC from the square of adjunctions
$$
\begin{matrix}
\left(\begin{CD} \O_{S_A}\backslash(A\FAlg(Y_0))\\  @VV{ f_0^{\sharp}}V \\ \O_{S_A}\backslash(A\FAlg(X_0)) \end{CD}\right)&\begin{smallmatrix} \xra{\quad\quad\quad\quad   }\\ \top \\ \xla[\O_{\fS}\ten_{\L}\Symm_{A}]{} \end{smallmatrix}&\left(\begin{CD} A\FMod(Y_0)\\@VV{ f_0^{\sharp}}V \\ A\FMod(X_0)\end{CD}\right)\\
\\
	u_*\uparrow \vdash\downarrow u^{-1}&	&u_*\uparrow \vdash\downarrow u^{-1}\\
\\
\left(\begin{CD}\O_{S_A}\backslash(A\FAlg(Y_0'))\\ @VV{ (f_0')^{\sharp}}V\\ \O_{S_A}\backslash(A\FAlg(X_0'))\end{CD}\right)&\begin{smallmatrix} \xra{\quad\quad\quad\quad   }\\ \top \\ \xla[\O_{\fS}\ten_{\L}\Symm_{A}]{} \end{smallmatrix}& \left(\begin{CD}A\FMod(Y_0')\\  @VV{ (f_0')^{\sharp}}V \\  A\FMod(X_0')\end{CD}\right),
\end{matrix}
$$
where for categories $\cA,\cB$, with a functor $F:\cA \to \cB$, $(\cA \xra{F} \cB)$ is the category whose objects are triples
$$
(A,B, FA \xra{h} B),
$$
and whose morphisms are pairs $f:A \to A'$, $g:B \to B'$ such that the following square commutes:
$$
\begin{CD} FA  @>{h}>>  B\\
 @V{Ff}VV   @VV{g}V\\
FA' @>{h'}>>  B'. \end{CD}
$$

\subsection{Deformation of a group scheme}\label{gpscheme}

Given a group scheme $G_0/k$, we consider deformations 
$$
\begin{CD}
G_0				@>>> 	G_A		\\
@VVV			 	@VV{\text{flat}}V\\
\Spec k 			@>>>	 \Spec A, 
\end{CD}
$$
where $G_A$ is a group scheme over $A$ with $G_A \ten_{A}k=G_0$.

Since the topological space $|G_0|$ does not deform, we need only consider deformations of  $\O_{G_0}$  as a sheaf of Hopf algebras. 
The structures are:
$$
\Sigma^+=\{\text{Algebra}\},\quad \Sigma^-=\{G_0\text{-Sheaf},\,\, \text{Co-Algebra}\},
$$
over the category of sheaves of flat $A$-modules on $G_0'$.

This gives the following commutative diagram of monadic and comonadic adjunctions:
$$
\xymatrix@C=12ex@R=8ex{
A\mathrm{-FHopfAlg}(G_0) \ar@<1ex>[r]_-{\top} \ar@<-1ex>[d]_{u^{-1}}
& \ar@<1ex>[l]^-{\Symm_{A}}	A\mathrm{-FCoAlg}(G_0)\ar@<-1ex>[d]_{u^{-1}}
\\
\ar@<-1ex>[u]_{Q\circ u_*}^{\dashv} A\FAlg(G_0') \ar@<1ex>[r]_-{\top}
&\ar@<-1ex>[u]_{Q \circ u_*}^{\dashv} \ar@<1ex>[l]^-{\Symm_{A}} A\FMod(G_0'),}
$$
where 
$A\mathrm{-FHopfAlg}$ is the category of  flat Hopf algebras over $A$, $A\mathrm{-FCoAlg}$ is the category of  flat co-associative co-algebras with co-unit and co-inverse over $A$, and
$Q$ is the free co-algebra functor of Section \ref{coalg}.  
We thus obtain the SDC
$$
E^n=\Hom_{\L}((\Symm_{\L})^n \sM, (u_*\circ Q\circ u^{-1})^n \sM)_{u^{-1}(\alpha^n \circ \vareps^n)},
$$
where $\sM$ is the lift (unique up to isomorphism) of the sheaf $u^{-1}\O_{G_0}$ of  $k$-vector spaces to a flat $\mu$-adic sheaf on $G'$, 
with $\alpha^n, \vareps^n$  the canonical maps 
$$
\vareps^n:(\Symm_{k})^n \O_{G_0} \to \O_{G_0},\quad
\alpha^n: \O_{G_0}  \to (u_*\circ Q\circ u^{-1})^n \O_{G_0}
$$
associated to $\O_{G_0}$.   

\subsubsection{Deformation of a sub-(group scheme)}
The problem is now, given  a flat $\mu$-adic formal group scheme $\fH$ over $\L$,  to deform a  sub-(group scheme) $G_0 \to H_0$ as a sub-(group scheme) of $H_A := \fH \ten_{\L}A$, flat over $A$. In fact, we will not need to assume that $G_0 \to H_0$ is a subscheme --- any morphism will do. The picture is:
$$
\xymatrix{
G_0 \ar[dr] \ar[r] \ar[rrd] &H_0 \ar[d] \ar[rrd]\\
&\Spec k \ar[rrd] &G_A \ar[dr]|{\text{flat}} \ar[r] & H_A \ar[d]^{\text{flat}}\\
&&&\Spec A,
}
$$
all morphisms as group schemes over $\L$.

The structures are:
$$
\Sigma^+=\{\text{Algebra},\,\,\O_{H_A}\text{-Module}\},\quad \Sigma^-=\{G_0\text{-Sheaf},\,\, \text{Coalgebra}\},
$$
so the diagram
$$
\xymatrix@C=12ex@R=8ex{
\O_{H_A}\backslash (A\mathrm{-FHopfAlg}(G_0)) \ar@<1ex>[r]_-{\top}\ar@<-1ex>[d]_{u^{-1}}
&\ar@<1ex>[l]^-{\O_{\fH}\ten_{\L}\Symm_{A}}  A\mathrm{-FCoAlg}(G_0)\ar@<-1ex>[d]_{u^{-1}}
\\
\ar@<-1ex>[u]_{Q\circ u_*}^{\dashv} \O_{H_A}\backslash (A\FAlg(G_0')) \ar@<1ex>[r]_-{\top} 
& \ar@<1ex>[l]^-{\O_{\fH}\ten_{\L}\Symm_{A}}  
\ar@<-1ex>[u]_{Q \circ u_*}^{\dashv}  A\FMod(G_0')}
$$
gives us our SDC, 
$$
E^n=\Hom_{\L}((\O_{\fH}\ten_{\L}\Symm_{\L})^n \sM, (u_*\circ Q\circ u^{-1})^n \sM)_{u^{-1}(\alpha^n \circ \vareps^n)}.
$$

\subsection{$G$-invariant deformations}
Given  a flat $\mu$-adic formal group scheme $\fG$ over $\Spf \L$, a flat $\mu$-adic formal scheme $\fT/\L$, and a flat scheme $X_0/k$, with a morphism $X_0 \to T_0$, and compatible group actions
$$
\fT\hat{\by}_{\L}\fG \to \fT,  \quad r:X_0 \by_k G_0 \to X_0,
$$
where $G_0=\fG \ten_{\L}k$, we seek deformations
$$
\xymatrix{
X_0 \ar[dr] \ar[r] \ar[rrd] &T_0 \ar[d] \ar[rrd]\\
&\Spec k \ar[rrd] &X_A \ar[dr]|{\text{flat}} \ar[r] & T_A \ar[d]^{\text{flat}}\\
&&&\Spec A,
}
$$
where $T_A=\fT \ten_{\L}A$ and $G_A=\fG \ten_{\L}A$, together with a compatible group action $X_A \by_A G_A \to X_A$. Since the underlying topological space does not deform, we need only consider $G_A$-linearised sheaves $\sF$ of algebras, i.e sheaves   with co-associative maps \mbox{$\sF \to r_*(\sF\ten_A \O_{G_A})$,}    deforming $\O_{X_0}$. 

The structures are:
$$
\Sigma^+=\{\text{Algebra},\,\,\O_{T_A}\text{-Module}\},\quad \Sigma^-=\{X_0\text{-Sheaf},\,\, G_A\text{-linear} \},
$$
where a $G_A$-linearised module is a  module $M$ equipped with a co-associative morphism
$$
M \to M \ten \O_{G_A}.
$$
 
The relevant diagram is then:
$$
\xymatrix@C=12ex@R=8ex{
\O_{T_A}\backslash((A\FAlg(X_0))^{G_A}) \ar@<1ex>[r]_-{\top} \ar@<-1ex>[d]_{u^{-1}}
&\ar@<1ex>[l]^-{\O_{\fT}\ten_{\L}\Symm_{A}}    (A\FMod(X_0))^{G_A} \ar@<-1ex>[d]_{u^{-1}}
\\
\ar@<-1ex>[u]_{R \circ u_*}^{\dashv} 	\O_{T_A}\backslash(A\FAlg(X_0')) \ar@<1ex>[r]_-{\top} 
&\ar@<-1ex>[u]_{R \circ u_*}^{\dashv} \ar@<1ex>[l]^-{\O_{\fT}\ten_{\L}\Symm_{A}} A\FMod(X_0'),}
$$
where $(A\FMod)^{G_A}$ denotes $G_A$-linearised module, and $R$ is the co-free functor
$$
R(\sF)= r_{*}(\sF \ten_{\L} \O_{\fG}),
$$ 
whose $\fG$-action is given by:
$$
r_{*}(\sF \ten_{\L} \O_{\fG}) \xra{r_*(\id \ten m^{\sharp})} r_*(\sF \ten m_*(\O_{\fG} \ten \O_{\fG}))=r_*(r_{*}(\sF \ten_{\L} \O_{\fG})\ten \O_{\fG}),
$$
as required.
The unit of the adjunction $V \dashv R$, for forgetful functor $V$, is 
$$
\rho: \sF \to  RV(\sF)=r_{*}(\sF \ten_{\L} \O_{\fG}),
$$
where $\rho$ is the structure map of $\sF$ as a $\fG$-linearised sheaf. The co-unit is 
$$
r_*(\id \ten e^{\sharp}): VR(\sF)=r_{*}(\sF \ten_{\L} \O_{\fG}) \to \sF.
$$

This gives the SDC
$$
E^n=\Hom_{\L}((\O_{\fT}\ten_{\L}\Symm_{\L})^n \sM, (u_*\circ R\circ u^{-1})^n \sM)_{u^{-1}(\alpha^n \circ \vareps^n)},
$$
where $\sM$ is the lift (unique up to isomorphism) of the sheaf $u^{-1}\O_{X_0}$ of  $k$-vector spaces to a flat $\mu$-adic sheaf on $X'$.

\section{Representations of the Algebraic Fundamental Group}\label{rep}

Fix a connected scheme $X$, a geometric point $\bar{x} \to X$, and a  ring $\L$ as in Section \ref{artin}, with finite residue field $k$. Throughout, $A$ will denote a ring in $\C_{\L}$. Denote $\pi_1(X,\bar{x})$ by $\Gamma$.

\subsection{Representations to $\GL_n$}\label{gln}

It is well known (e.g. \cite{Mi} Theorem 5.3), that there is an equivalence of categories between finite $\pi_1(X,\bar{x})$-sets, and locally constant \'etale sheaves on $X$, with finite stalks. There is thus an equivalence of categories between representations 
$$
\rho: \pi_1(X,\bar{x}) \to \GL_n(A),
$$ and locally constant \'etale sheaves $\vv$ on $X$ with stalks $A^n$.  Given $\rho$, the corresponding sheaf $\vv_{\rho}$ is given by 
$$
\vv_{\rho}=A^n\by_{\Gamma, \rho}\F_Y=(\pi_* A^n)^{\Gamma,\rho},
$$ 
where $\F_Y$ is the sheaf on $X$ whose espace \'etal\'e is $Y$, and $\rho|_{\pi_1(Y,\bar{y})}=1$, which must happen for some  Galois cover $Y\xra{\pi} X$, $A$ being finite.

If we now fix a representation $\rho_0:\Gamma \to \GL_n(k)$, then we have an equivalence of categories (or, rather, groupoids) 
$$
\r_A(\rho_0) \leftrightarrow \mathfrak{V}_A(\vv_0),
$$ 
where $\r_A(\rho_0)$ has objects 
$$\rho \equiv \rho_0 \mod \m_A, 
$$
while $\mathfrak{V}_A(\vv_0)$ has objects 
$$\vv \equiv \vv_0 \mod \m_A.
$$ 
In both cases morphisms must be the identity mod $ \m_A$. Explicitly, a morphism in $\r_A(\rho_0)$ is an element $g \in \ker(\GL_n(A) \to \GL_n(k))$, acting by conjugation.

Hence  this deformation problem is just the same as deforming sheaves of modules, as in Section \ref{etshf}, giving the SDC
$$
E^n=\Hom_{\L}(\sV, (u^*u_*)^n(\sV))_{u^*\alpha^n},
$$
where $\sV$ is the lift (unique up to isomorphism) of the sheaf $u^{-1}\vv_0$ of  $k$-vector spaces to a flat $\mu$-adic sheaf on $X'$. The cohomotopy groups are, of course, isomorphic to
$$
\H^*(X,\hom_k(\vv_0,\vv_0)).
$$ 
Note that 
$$
\hom_k(\vv_0,\vv_0)=(\pi_*\mathfrak{gl}_n(k))^{\Gamma,\rho_0},
$$
so the cohomotopy groups are
$$
\H^i(X,(\pi_*\mathfrak{gl}_n(k))^{\Gamma,\rho_0}),
$$
where $\Gamma$ acts via the adjoint action.

\subsection{Representations  to a smooth Noetherian group scheme}\label{smoothgroup}

Fix a system  $\{G_n/\Spec \L_n\}_{n \in \N}$ of smooth Noetherian  group schemes, and let $G$ be the $\mu$-adic formal group scheme $\varinjlim G_n$.   Fix a representation $\rho_0: \Gamma \to G(k)$.

\begin{definition} A principal $G(A)$-sheaf on $X$ is an \'etale sheaf $\Bu$ on $X$, together with an action $ \G\by \Bu  \to \Bu$, where $\G$ denotes the constant sheaf $G(A)$, such that 
$$
\G\by\Bu \cong \Bu\by\Bu,\quad (g,b)\mapsto (g\cdot b,b),
$$
and the stalks $\Bu_{\bar{x}}$ are non-empty.

Given a map $A \to B$ in $\C_{\L}$, there is a  functorial map from $G(A)$-sheaves to $G(B)$-sheaves, given by
$$
\Bu \mapsto G(B)\by^{G(A)}\Bu,
$$ 
i.e. $(gh,b)\sim (h,gb)$, for $g \in G(A)$, $h \in G(B)$ and $b \in B$.
\end{definition}

\begin{definition} Given a  representation $\rho: \Gamma \to G(A)$, let
$$
\Bu_{\rho}:=\F_Y\by^{\Gamma, \rho}G(A)=(\pi_*G(A))^{\Gamma,\rho},
$$ 
where $\F_Y$ is the sheaf on $X$ whose espace \'etal\'e is $Y$, and $\rho|_{\pi_1(Y,\bar{y})}=1$, which must happen for some  Galois cover $Y \xra{\pi}X$, $G(A)$ being finite (since $G$ is Noetherian and $A$ finite).

More precisely, 
$$
\F_Y(U \xra{u} X)=\left\{
\vcenter{\xymatrix@=3ex{
U \ar[r]^{\tilde{u}} \ar[dr]_u  &Y \ar[d] \\
&X
}}\quad \text{a lift}\right\},
$$
and the equivalence relation forming $\Bu_{\rho}$ is $(\gamma \circ u,g)\sim (u, g\rho(\gamma))$, or equivalently \mbox{$\eta \sim\gamma^*\eta\rho(\gamma) $.}
$\Bu_0:=\{\Bu_{\rho_0|_{\L_n}}\}_{n\in \N}$.
\end{definition}

\begin{lemma} There is an equivalence of categories
$$
\r_A(\rho_0) \leftrightarrow \mathfrak{B}_A(\Bu_0),
$$ 
where $\r_A(\rho_0)$ has as  objects  representations 
$$
\rho \equiv \rho_0 \mod \m_A,
$$ 
while $\mathfrak{B}_A(\Bu_0)$ has as objects principal $G(A)$-sheaves  
$$
\Bu \equiv \Bu_0 \mod \m_A
$$ on $X$. In both cases morphisms must be the identity mod $ \m_A$. 

\begin{proof}
We need to show that the functor $\Phi:\rho \mapsto \Bu_{\rho}$ is full, faithful and essentially surjective.
\begin{enumerate}
\item $\Phi$ is full and faithful:

Given $\theta: \Bu_{\sigma} \xra{\sim} \Bu_{\tau}$, take some finite Galois cover $Y \xra{\pi} X$ such that $\pi^{-1}\Bu_{\sigma}$ and $\pi^{-1}\Bu_{\tau}$ are both isomorphic to the trivial sheaf  $G(A)$. Thus $\pi^{-1}(\theta)$ acts as a unique element $\alpha$ of $G(A)$ with  $\alpha \equiv 1 \mod \m_A$, so $\theta = \Phi(\alpha)$.

\item $\Phi$ is essentially surjective:

Given $\Bu \in \mathfrak{B}_A(\Bu_0)$, $\Bu$ is locally constant with finite stalks (since  $G(A)$ is finite), so there is a Galois cover  $Y \xrightarrow{\pi} X$ such that $\pi^{-1}\Bu$ is constant.  Choose a geometric point  $\bar{y} \to Y$, above $\bar{x}$. $\Gamma$ acts on $(\pi^{-1}\Bu)_{\bar{y}}$ via:
$$ (\pi^{-1}\Bu)_{\bar{y}} \cong (\pi^{-1}\Bu)_{\gamma\circ \bar {y}} \cong \Bu_{\bar{x}} \cong (\pi^{-1}\Bu)_{\bar{y}},$$
the first isomorphism being given by the constancy of $\pi^{-1}\Bu$, and the other two canonical. This is an isomorphism of $G(A)$-sets , so gives an element $\rho \in \r_A(\rho_0) $. It follows that $\Phi(\rho) \cong \Bu$.
\end{enumerate}
\end{proof}
\end{lemma}

Observe that deformations of $\Bu_0$ as a principal $G$-sheaf are the same as its deformations as a faithful $G$-sheaf.
We have the functorial comonadic adjunction:
$$
\xymatrix{
\text{Faithful }G(A)\text{-Sheaves}(X) \ar@<-1ex>[d]_{u^{-1}}\\
\text{Faithful }G(A)\text{-Sheaves}(X') \ar@<-1ex>[u]_{u_*}^{\dashv}
}
$$
using the canonical injection $G(A) \to u_* G(A)$ to make $u_*\sB$ a faithful $G(A)$-sheaf.
Observe that the second category has uniformly trivial deformation theory, $G$ being smooth, so that deformations are described by the SDC
$$
E^n=\Hom(\B,(u^*u_*)^n\B)^{G(\L)}_{u^*\alpha^n},
$$
where $\sB$ is the lift (unique up to isomorphism) of the sheaf $u^{-1}\Bu_0$ of  $G(k)$-sheaves to a sheaf of continuous (in the $\mu$-adic topology) $G(\L)$-sheaves on $X'$. The cohomotopy groups are isomorphic to 
$$
\H^*(X,\bu_0),
$$ 
where $\bu_0$ is the tangent space of $\aut(\Bu_0)$ at the identity. Equivalently,
$$
\bu_0= \g \by^{G(k)}\Bu_0=(\pi_* \g)^{\Gamma,\rho_0},
$$
for $\g$ the tangent space of $G(k)$ at the identity, and $G(k)$ acting on $\g$ via the adjoint action,
so the cohomotopy groups are
$$
\H^*(X,(\pi_*\g)^{\Gamma,\rho_0}).
$$

\section{Equivalence between SDCs and DGLAs in Characteristic Zero}\label{sdcdgla}

Throughout this section, we will assume that $\L=k$, a field of characteristic zero. We will also assume that all DGLAs are in non-negative degrees, i.e. $L^{<0}=0$.

\begin{remark} Throughout this section, we will encounter both chain complexes and simplicial complexes. In order to avoid confusion, chain complexes will be denoted by a bullet ($C_{\bullet}$), while simplicial complexes will be denoted by an asterisk ($C_*$). Asterisks will also be used to denote complexes without differentials, so that, for a chain complex $C_{\bullet}$, we may write $C_{\bullet}=(C_*,d)$.  Homology groups of a chain complex $C_{\bullet}$ will be denoted $\H_i(C)$, while homotopy groups of a simplicial complex $C_*$ will be denoted $\pi_i(C)$. The dual conventions will be used for cochain and cosimplicial complexes.
\end{remark}

\begin{enumerate}
\item Consider, for simplicity, an SDC $E^*$ for which each tangent space $\mathfrak{e}^n:=\CC^n(E)$ is finite dimensional. Thus each $E^n$ is pro-representable, by a  smooth complete local  Noetherian $k$-algebra $Q_n$, say. Since $\pd^0_{\omega_0}$ and $\pd^{n+1}_{\omega_0}$ are defined on $\C_k$, $E^*$ is a cosimplicial complex, so  $Q_*$ has the structure of a simplicial complex of algebras. The product on $E^*$ gives us a comultiplication 
$$
*^{\sharp}_{m,n}:Q_{m+n} \to Q_m \hat{\ten} Q_n,
$$
satisfying various axioms. $\mathfrak{e}^*=\Der_k(Q_*,k)$, so $\pi^i(\mathfrak{e}^*)=\pi^i(\Der_k(Q_*,k))$.

\item To give a graded Lie algebra $L^*$ (in non-negative degrees) over $k$ is the same as giving the smooth homogeneous functor
\begin{eqnarray*}
\exp(L): \C_k^{\N_0} &\to& \Grp\\
		A_* &\mapsto& \exp(\bigoplus_n L^n \ten (\m_A)_n),
\end{eqnarray*}
where $\C_k^{\N_0}$ is the category of nilpotent local Noetherian $\N_0$-graded (super-commutative) $k$-algebras with residue field $k$ (concentrated in degree $0$). Observe that $\bigoplus_n L^n \ten (\m_A)_n$ is a Lie algebra, so the formula makes sense.

The category $\C_k^{\N_0}$ behaves similarly to $\C_k$, in that the analogue of Theorem \ref{Sch} holds; in particular, if $L^*$ is finite dimensional, then $\exp(L)$ can be pro-represented by a smooth complete local Noetherian graded $k$-algebra $R_*$, the group structure making $R_*$ into a graded Hopf algebra:
$$
*^{\sharp}_{m,n}:R_{m+n} \to R_m \hat{\ten} R_n.
$$

If we take a finite dimensional DGLA $L^{\bullet}$, the the differential gives rise to a differential on $R_*$, giving it the structure of a DG Hopf algebra $R_{\bullet}$. Since $L^{\bullet}= \Der_k(R_{\bullet},k)$, we have $\H^i(L^{\bullet})= \H^i(\Der_k(R_{\bullet},k)$.
\end{enumerate}

Hence, given the equivalence, described in \cite{QRat}, between DG and simplicial algebras in characteristic zero, it should be possible to obtain DGLAs from SDCs and \emph{vice versa}. Moreover, given the similarity between our cohomology groups and Andr\'e-Quillen cohomology, and the result in \cite{Q} that, in characteristic 0, Andr\'e-Quillen cohomology can be computed by either simplicial or DG resolutions, it seems that an SDC should have the same cohomology as the corresponding DGLA. It is then plausible that they should have the same deformation groupoid. The purpose of this section is to make this reasoning rigorous.

\subsection{SDCs}\label{sdcrep}

Fix an SDC $E^*$. Although each $E^n$ will only be pro-representable if $\dim \mathfrak{e}^n < \infty$, in general we may make use of the following lemma:

\begin{lemma}\label{proart}
Given a  smooth homogeneous functor $F:\C_{\L} \to \Set$, let $\{t_i:i \in I\}$ be a basis for $t_F$. Then $F$ is isomorphic to $\Hom_{\text{cts.}}(\cR, -)$, where $\cR$ is the pro-Artinian completion of $\L[T_i:i \in I]$, localised at $(T_i:i \in I)$, for formal symbols $T_i$. Note that, if $t_F$ is finite dimensional, this is just the expected ring $\L[[T_1,\ldots,T_m]]$. 

\begin{proof} This is very similar to the proof of Theorem \ref{Sch}. A full proof in the analogous case of nilpotent Lie algebras (rather than Artinian rings)  is given in \cite{higgs} Theorem \ref{higgs-nSch}.
\end{proof}
\end{lemma}

Thus an SDC $E^*$ is equivalent to a system of smooth local complete pro-Artinian rings $\cQ_*$, with the dual structures to those described in Section \ref{sdcdef}. In particular, since we are only considering the case $\L=k$, we have canonical $\pd^0_{\omega_0},\pd^{n+1}_{\omega_0}$. Explicitly:

\begin{lemma}\label{*rules}
$\cQ_*$ is a simplicial complex of smooth local complete pro-Artinian $k$-algebras (with residue field $k$), together with a co-associative  comultiplication
$$
\rho_{m,n}:=*^{\sharp}_{m,n}:\cQ_{m+n} \to \cQ_m \hat{\ten} \cQ_n,
$$
such that $\cQ_0 \to k$ is the co-unit,  satisfying:
\begin{enumerate}
\item $ (\pd_i\ten \id) \circ \rho_{m+1,n} = \rho_{m,n} \circ \pd_i$ for $i \le m$.
\item $(\id \ten \pd_i) \circ \rho_{m,n+1} = \rho_{m,n} \circ \pd_{m+i}$ for $i \ge 1$.
\item $ (\pd_{m+1}\ten \id) \circ \rho_{m+1,n} = (\id \ten \pd_0) \circ \rho_{m,n+1}$.
\item $ (\sigma_i\ten \id) \circ \rho_{m-1,n} = \rho_{m,n} \circ \sigma_i$.
\item $(\id \ten \sigma_i) \circ \rho_{m,n-1} = \rho_{m,n} \circ \sigma_{m+i}$.
\end{enumerate}
\end{lemma}
Of course, given a system $\cQ_*$ satisfying these conditions, we can clearly recover an SDC as $\Hom_{\text{cts.}}(\cQ_*, -)$, so that we may regard this as an equivalent definition of an SDC.

\begin{remark}\label{piq} We say that a simplicial complex $A_*$ of local $\L$-algebras is Artinian if it satisfies DCC on simplicial ideals. This is equivalent to saying that the algebras $A_n$ are Artinian, and the  normalised complex $N(A_*)$ is bounded.
It is easy to see that 
$
\cQ_* = \Lim Q_*,
$
the limit being over quotient \emph{complexes} $\cQ_* \to Q_*$  Artinian. Moreover, $\cQ_*$ is really just a shorthand for this inverse system. Thus, for any functor $F$, by $F(\cQ_*)$ we mean $\Lim F(Q_*)$. In particular, $\pi_*(\cQ_*):=\Lim \pi_*(Q_*)$. We will apply the same convention to the limits of DG complexes in the next section.
\end{remark}

\subsection{DGLAs}\label{dglarep}

Define $\C_k^{\N_0}$ to be  the category of nilpotent local Noetherian graded (super-commutative) $k$-algebras (in non-negative degrees) with residue field $k$ (concentrated in degree $0$).

A GLA $L^*$ over $k$ (in non-negative degrees) gives rise to a smooth homogeneous functor
\begin{eqnarray*}
\exp(L): \C_k^{\N_0} &\to& \Grp\\
		A_* &\mapsto& \exp(\bigoplus_n L^n \ten (\m_A)_n).
\end{eqnarray*}

Given a functor $F:\C_k^{\N_0} \to \Set$, define $t_F^i:=F(k[\eps_i])$, where $\eps_i$ is a variable in degree $i$, with $\eps_i^2=0$.

\begin{lemma}
Given a  smooth homogeneous functor $F:\C_k^{\N_0} \to \Set$, let $\{t_i:i \in I\}$ be a homogeneous  basis for $t_F$. Then $F$ is isomorphic to $\Hom_{\text{cts.}}(\cR, -)$, where $\cR$ is the pro-Artinian completion of the ring $\L[T_i:i \in I]_{\gd}$, localised at $(T_i:i \in I)$, for formal symbols $T_i \in \deg(t_i)$. Note that, if $t_F$ is finite dimensional, this is just the expected ring $\L[[T_1,\ldots,T_m]]_{\gd}$. 
\begin{proof}
The analogue of Theorem \ref{Sch} holds for graded rings, and this is just a generalisation of Lemma \ref{proart}.
\end{proof}
\end{lemma}

Let $\exp(L)$ be ``represented'' by the DG algebra $\cT_*$. The group structure of $\exp(L)$ corresponds to a graded Hopf algebra map
$$
\rho: \cT_* \to \Tot(\cT_* \hat{\ten} \cT_*),
$$
where 
$$
\Tot(\cT_* \hat{\ten} \cT_*)_n = \bigoplus_{r+s=n} \cT_r \hat{\ten} \cT_s,
$$ 
and the multiplication follows the usual (graded) convention: 
$$
(u \ten v ) (x \ten y) = (-1)^{\bar{v}\bar{x}} ux \ten vy,
$$
 for homogeneous elements. We can, of course, recover $L^*$ from $\cT_*$ by
$$
L^*=\Der_{\text{cts.}}(\cT_*,k).
$$
\begin{lemma}\label{expolie}
Conversely, given a  ``representable'' functor $G:\C_k^{\N_0} \to \Grp$, with tangent space $L^*$, there is a canonical isomorphism
$$
\exp(L) \cong G.
$$
In particular, this implies that $G$ must be smooth.
\begin{proof}
Let $G$ be ``represented'' by $\cT_*$. Then 
$$
G(A_*)=\Hom_{k\Alg}(\cT_*,A_*),\quad \quad  \bigoplus_n L^n \ten (\m_A)_n = \Der_k(\cT_*,A_*).
$$
We can embed both of these into $\Hom_k(\cT_*,A_*)$, on which we define the associative product $f*g= (f \ten g)\circ \rho$.
Now the maps
\begin{eqnarray*}
G(A_*) &\cong& \bigoplus_n L^n \ten (\m_A)_n\\
g      &\mapsto& \sum_{n \ge 1}(-1)^{n-1} \frac{(g-e)^{*n}}{n}\\
\sum_{n \ge 0} \frac{l^{*n}}{n!} &\mapsfrom& l,
\end{eqnarray*}
where $e$ is the identity map $\cT_* \to k$, give us our isomorphism. 
\end{proof}
\end{lemma}

Given a DGLA $L^{\bullet}$, it remains only to understand what the differential $d$ on $L^{\bullet}$ gives rise to on $\cT_*$.

To give $d$ is the same as giving the map 
$$
\xymatrix{
 L^* \ar[r]^-{\id +d\eps_1} \ar[dr]_{\id}	& L^*[\eps_1] \ar[d]\\
&L^*,
}
$$
and $\exp(L^*[\eps_1])$ is ``represented'' by $\cT_*[\eps_1]$. Thus $d$ corresponds to a derivation \mbox{$d^{\sharp}:\cT_{n+1} \to \cT_{n}$.} Since $d^2=0$, we get $(d^{\sharp})^2=0$. The condition that $d$ be a Lie algebra derivation corresponds to the statement
$$
\rho \circ d^{\sharp}= (\id \ten d^{\sharp} \pm d^{\sharp} \ten \id) \circ \rho,\text{ where }
$$
$$
(\pm d^{\sharp} \ten \id) (u \ten v) = (-1)^{\bar{v}} d^{\sharp}(u)\cdot v.
$$ 

In summary, a DGLA $L^{\bullet}$ is equivalent to a  smooth local complete pro-Artinian DG  $k$-algebra $\cT_{\bullet}$, with a co-associative chain map
$$
\rho: \cT_{\bullet} \to \Tot(\cT_{\bullet} \hat{\ten} \cT_{\bullet}),
$$
for which $\cT_{\bullet} \to k$ is the co-unit.

\subsection{SDC $\leadsto$ DGLA}\label{sdctodgla}

Given $\cQ_*$, as in Section \ref{sdcrep}, we  form a DG Hopf algebra $C(\cQ)_{\bullet}$, with
$
C(\cQ)_n= \cQ_n 
$
as a graded module. We let 
$$
d = \sum_{i=0}^{n+1} (-1)^i \pd_i:C(\cQ)_{n+1} \to C(\cQ)_n,
$$
 and give it an algebra structure via the Eilenberg-Zilber shuffle product: for $x \in C(\cQ)_m$ and $y \in C(\cQ)_n$, let
$$
x \nabla y = \sum_{(\mu,\nu) \in \mathrm{Sh}(m,n)} (-1)^{(\mu,\nu)} (\sigma_{\nu_n}\ldots \sigma_{\nu_1}x) \cdot (\sigma_{\mu_m}\ldots \sigma_{\mu_1}y),
$$
where $\mathrm{Sh}(m,n)$ is the set of $(m,n)$ shuffle permutations, i.e. permutations 
$$
(\mu_1,\ldots,\mu_m,\nu_1,\ldots,\nu_n)\quad\text{ of }\quad\{0,\ldots,m+n-1\}
$$
 such that 
$$
\mu_1<\ldots <\mu_m \text{ and }\nu_1<\ldots<\nu_n,
$$
 and $(-1)^{(\mu,\nu)}$ denotes the sign of the permutation $(\mu,\nu)$. As explained in \cite{W} or \cite{QRat}, the resulting product is supercommutative 
$$
x \nabla y = (-1)^{mn} y \nabla x
$$ 
and associative, and $d$ is a graded derivation with respect to $\nabla$, i.e.
$$
d(x \nabla y) = (d x) \nabla y +(-1)^m x \nabla (dy).
$$
\begin{lemma}\label{*shuffle}
$C(\cQ)_*$ is a Hopf algebra under the map $\rho$.
\begin{proof}
 Given $x \in C(\cQ)_m$,  $y \in C(\cQ)_n$, and any $p+q=m+n$, we have:
\begin{eqnarray*}
\rho_{p,q}(x \nabla y)&=& \rho_{p,q}(\sum_{(\mu,\nu) } (-1)^{(\mu,\nu)} \sigma_{\nu_n}\ldots \sigma_{\nu_1}x \cdot \sigma_{\mu_m}\ldots \sigma_{\mu_1}y)\\
&=&\sum_{(\mu,\nu) } (-1)^{(\mu,\nu)}\rho_{p,q}(\sigma_{\nu_n}\ldots \sigma_{\nu_1}x)\cdot \rho_{p,q}(\sigma_{\mu_m}\ldots \sigma_{\mu_1}y)\\
&=& \sum
\pm (\sigma_{\nu'_{n'}}\ldots \sigma_{\nu'_1} \ten \sigma_{\nu''_{n''}}\ldots \sigma_{\nu''_1})\circ \rho_{m',m''}(x) \cdot (\sigma_{\mu'_{m'}}\ldots \sigma_{\mu'_1} \ten \sigma_{\mu''_{m''}}\ldots \sigma_{\mu''_1})\circ \rho_{n',n''}(y),
\end{eqnarray*}
the final sum being over
\begin{eqnarray*}
&m'+m''=m,\quad n'+n''=n,\quad m'+n'=p,\quad m''+n''=q,\\ 
&(\mu',\nu') \in \mathrm{Sh}(m',n'),\quad (\mu'',\nu'') \in \mathrm{Sh}(m'',n''),
\end{eqnarray*}
with
$$
\pm=(-1)^{(\mu',\mu''+p,\nu', \nu'' +p)}.
$$

For this map to be a Hopf algebra map, we would need this expression to equal
$$
\sum_{\begin{smallmatrix}m'+m''=m,n'+n''=n,\\m'+n'=p,m''+n''=q\end{smallmatrix}} \rho_{m',m''}(x) \nabla \rho_{n',n''}(y),
$$
The $(p,q)$ component of $\rho(x)\nabla \rho(y)$.
Now,
\begin{eqnarray*}
&&\rho_{m',m''}(x) \nabla \rho_{n',n''}(y)=\\
 &&\sum_{\begin{smallmatrix}(\mu',\nu')\\(\mu'',\nu'')\end{smallmatrix}}
\pm (\sigma_{\nu'_{n'}}\ldots \sigma_{\nu'_1}\ten \sigma_{\nu''_{n''}}\ldots \sigma_{\nu''_1})\circ \rho_{m',m''}(x)  \cdot (\sigma_{\mu'_{m'}}\ldots \sigma_{\mu'_1} \ten \sigma_{\mu''_{m''}}\ldots \sigma_{\mu''_1})\circ \rho_{n',n''}(y),
\end{eqnarray*}
where 
$$
\pm=(-1)^{m''n'}(-1)^{(\mu',\nu')}(-1)^{(\mu'',\nu'')},
$$
the $(-1)^{m''n'}$ arising from the (required) supercommutativity of $C(\cQ)_*\ten C(\cQ)_*$. The result now follows from the observation that 
$$
(-1)^{(\mu',\mu''+p,\nu', \nu'' +p)}=(-1)^{m''n'}(-1)^{(\mu',\nu')}(-1)^{(\mu'',\nu'')}.
$$
\end{proof}
\end{lemma}

\begin{lemma}\label{*d}
$\rho$ is a chain map.
\begin{proof} For $x \in C(\cQ)_{m+1}$,
\begin{eqnarray*}
\rho(d x)&=&\sum_{p+q=m} \rho_{p,q}(\sum_{i=0}^{m+1} (-1)^i \pd_i x)\\
&=& \sum_{p+q=n}( \sum_{i=0}^p(-1)^i (\pd_i \ten \id)\circ \rho_{p+1,q}(x) + (-1)^p \sum_{j=1}^{q+1}(-1)^j (\id \ten \pd_j)\circ \rho_{p,q+1}(x))\\
&=&\sum_{p+q=m}( ((d-(-1)^{p+1}\pd_{p+1}) \ten \id) \circ \rho_{p+1,q}(x) +(-1)^p (\id \ten (d-\pd_0))\circ \rho_{p,q+1}(x))\\
&=&  \sum_{p+q=m}(d \ten \id)\circ \rho_{p+1,q}(x) -(-1)^{p+1}(\pd_{p+1} \ten \id)\circ\rho_{p+1,q}(x)\\&& +  \sum_{p+q=m}(-1)^p  (\id \ten d)\circ \rho_{p,q+1}(x) -(-1)^p (\id \ten \pd_0)\circ \rho_{p,q+1}(x)\\
&=& \sum_{p+q=m}(d \ten \id)\circ \rho_{p+1,q}(x) +   \sum_{p+q=m}(-1)^p  (\id \ten d)\circ \rho_{p,q+1}(x)\\
&=& d\rho(x).
\end{eqnarray*}
\end{proof}
\end{lemma}

There is now a problem with $C(\cQ_*)$ as a candidate to define a DGLA $L$ --- it is not local. In \cite{QRat}, rather than taking the unnormalised chain complex $C(\cQ_*)$, Quillen takes the normalised chain complex $N(\cQ_*)$. Recall that
$$
C(\cQ)_n=N(\cQ)_n \oplus D(\cQ)_n, \text{ where }
$$
$$
N(\cQ)_n = \bigcap_{i=1}^{n} \ker(\pd_i: \cQ_n \to \cQ_{n-1}),\quad \quad D(\cQ)_n= \sum_{i=1}^{n-1} \sigma_i(\cQ_{n-1}).
$$
Observe that $N(\cQ)_{\bullet}$ is local, as $\ker(\pd_1) \subset \m_{\cQ}$. It is pro-Artinian following Remark \ref{piq}, since the normalisation of an ideal is a $\nabla$-ideal, and $N(I)\nabla N(J)\subset N(IJ)$. However, we cannot give $N(\cQ)$ the Hopf algebra structure $\rho$: in general, for $x \in N(\cQ)_{m+n}$, $(\pd_m \ten \id)\circ \rho_{m,n}(x) \ne 0$. Instead, we take the complex 
$$
\bar{N}(\cQ):=C(\cQ)/D(\cQ).
$$ 
It follows immediately from the identities in Lemma \ref{*rules} that 
$$
\rho_{m,n}:\bar{N}(\cQ)_{m+n} \to \bar{N}(\cQ)_{m} \hat{\ten} \bar{N}(\cQ)_n
$$
is well defined. We need only check that $\nabla$ is well defined on $\bar{N}(\cQ)$, since then we will have $\bar{N}(\cQ)$ isomorphic to $N(\cQ)$ as an algebra.

\begin{lemma}
For $x \in D(\cQ)_m$ and $y \in C(\cQ)_n$, $x \nabla y \in D(\cQ)_{m+n}$.
\begin{proof}
Without loss of generality, $x = \sigma_r u$, for some $0 \le r < m$. We will show that each term in the sum of the shuffle product is in $D(\cQ)_{m+n}$. Consider
$$
\sigma_{\nu_n}\ldots \sigma_{\nu_1}x \cdot \sigma_{\mu_m}\ldots \sigma_{\mu_1}y,
$$
and let $ a := \max \{j: \nu_j -j < r\}$, which makes sense, since $\nu_j - j$ is a non-decreasing function. With the obvious convention if this set is empty, we have $0 \le a \le n$. Now, the simplicial identities give
$$
\sigma_{\nu_n}\ldots \sigma_{\nu_1}\sigma_r u = \sigma_{\nu_n}\ldots \sigma_{\nu_{a+1}} \sigma_{r+a} \sigma_{\nu_a}\ldots \sigma_{\nu_1} u,
$$
with 
$$
\nu_1<\ldots<\nu_a < r+a < \nu_{a+1}<\ldots<\nu_n.
$$
Now $r+a \le p-1 +q$, and $r+a$ does not equal any of the $\nu_j$, so $r+a=\mu_i$, for some $i$. This gives
\begin{eqnarray*}
(\sigma_{\nu_n}\ldots \sigma_{\nu_1}x) \cdot (\sigma_{\mu_m}\ldots \sigma_{\mu_1}y) &=&(\sigma_{\nu_n}\ldots \sigma_{\nu_{a+1}} \sigma_{\mu_i} \sigma_{\nu_a}\ldots \sigma_{\nu_1} u  )\cdot( \sigma_{\mu_m}\ldots \sigma_{\mu_1}y)\\
&=&( \sigma_{\mu_i}\sigma_{\nu_n-1}\ldots \sigma_{\nu_{a+1}-1} \sigma_{\nu_a}\ldots \sigma_{\nu_1} u  )\cdot(\sigma_{\mu_i}\sigma_{\mu_m-1}\ldots \sigma_{\mu_{i+1}-1}\sigma_{\mu_{i-1}}\ldots \sigma_{\mu_1} y)\\
&=&  \sigma_{\mu_i}((\sigma_{\nu_n-1}\ldots \sigma_{\nu_{a+1}-1} \sigma_{\nu_a}\ldots \sigma_{\nu_1} u  )\cdot(\sigma_{\mu_m-1}\ldots \sigma_{\mu_{i+1}-1}\sigma_{\mu_{i-1}}\ldots \sigma_{\mu_1} y)),
\end{eqnarray*}
which is in $D(\cQ)_{m+n}$, as required.
\end{proof}
\end{lemma}

\begin{definition} Given an SDC $E^*$, we define the DGLA $\cL(E)^{\bullet}$ by
$$
\cL(E)^{\bullet} = \Der_{k,\text{cts.}}(\bar{N}(\cQ)_{\bullet},k),
$$
where $\cQ_*$ ``represents'' $E^*$. 
\end{definition}

\begin{lemma}\label{hgtol} There are natural maps $\pi^i(\mathfrak{e}^*) \to \H^i(\cL(E)^{\bullet})$, where $\mathfrak{e}^*$ is the tangent space of $E^*$.
\begin{proof} 
We may compute $\pi^i$ using the conormalised cocomplex 
$$
N^n(K^*)= \bigcap_{i=0}^{n-1} \ker(\sigma^i: K^n \to K^{n-1}),\quad  d = \sum_{i=0}^n (-1)^i \pd^i.
$$
 Thus an element $\alpha \in N^i(\mathfrak{e}^*)$ is a derivation (with respect to the product on $\cQ$) from $\bar{N}_i(\cQ)$ to $k$ (by the definition of $\bar{N}$). A derivation to $k$ is just something which annihilates both $k$ and $\m_{\cQ}^2$, but $\m_{\cQ}\nabla \m_{\cQ} \subset \m_{\cQ}^2$, so $\alpha$ gives a derivation (with respect to $\nabla$) from  $\bar{N}_{\bullet}(\cQ)$ to $k$, of degree $i$, i.e. 
$$
\alpha  \in \cL(E)^i.
$$
 This correspondence preserves closed and exact forms, so we get 
$$
\pi^i(\mathfrak{e}^*) \to \H^i(\cL(E)^{\bullet}).
$$
\end{proof}
\end{lemma}

\begin{lemma}\label{gtol}
There is a natural map $\Def_E \to \Def_{\cL(E)}$, with the maps on cohomology compatible with the corresponding maps on tangent and obstruction spaces.
\begin{proof}
\begin{enumerate}
\item
Given $\omega \in \mc_E$, we may regard $\omega$ as a continuous ring homomorphism $\omega: \cQ_1 \to A$, satisfying $\omega*\omega=\pd^1(\omega)$. We also know that $\sigma^0(\omega)=e$. Set $\alpha=\omega -\omega_0$. This gives $\sigma^0(\alpha)=0$, therefore $\alpha:\bar{N}_1(\cQ) \to A$. In fact, $\alpha$ is a $\nabla$-derivation: for $x \in \cQ_0, y \in \cQ_1$,
\begin{eqnarray*}
\alpha(x\nabla y) &=& \alpha (\sigma_0(x)y)\\
&=& (\omega -\omega_0)(\sigma_0(x)y)\\
&=& \sigma^0(\omega)(x)\bar{y} - \sigma^0(\omega_0)(x)\bar{y}  + \overline{\sigma_0(x)} \omega(y)- \overline{\sigma_0(x)} \omega_0(y)\\
&=& \bar{x}\alpha(y)\\
&=& \bar{x}\alpha(y) + \alpha(x)\bar{y}.
\end{eqnarray*}
Moreover, the Maurer-Cartan equation for $\omega$ becomes
\begin{eqnarray*}
(\omega_0 +\alpha)*(\omega_0 +\alpha) &=& \pd^1(\omega_0+\alpha)\\
\omega_0^2 + \pd^0(\alpha) +\pd^2(\alpha) +\alpha*\alpha &=& \omega_0^2 +\pd^1(\alpha)\\
\pd^0(\alpha)-\pd^1(\alpha) +\pd^2(\alpha)+\alpha*\alpha &=& 0\\
d\alpha + \half [\alpha,\alpha]&=&0.
\end{eqnarray*} 
Therefore $\alpha  \in \mc_{\cL(E)}$.

\item Given $g \in E^0(A)$, $g$ is an algebra homomorphism $g:\cQ_0 \to A$. Since $\nabla$ agrees with the usual product on $\cQ_0$, and $\bar{N}_0(\cQ)=\cQ_0$, $g$ is a $\nabla$-algebra homomorphism $g:\bar{N}_0(\cQ) \to A$. Therefore 
$$
g  \in \exp(\cL(E)\ten \m_A).
$$
 It is easy to see that the adjoint action on $E^*$ corresponds to the gauge action on $\cL(E)^{\bullet}$.
\end{enumerate}

We therefore have the map $\Def_E \to \Def_{\cL(E)}$. The results on tangent and obstruction spaces follow immediately.
\end{proof}
\end{lemma}

We now need only  prove that the maps $\pi^i(\mathfrak{e}^*) \to \H^i(\cL(E)^{\bullet})$ are isomorphisms, from which it will follow that the morphism of groupoids is an equivalence. This will be proved in Theorem \ref{gtol2}.

\subsection{DGLA $\leadsto$ SDC}

\begin{definition}
Given a chain complex $V_{\bullet}$ of vector spaces, the left adjoint $N^{-1}$ of the normalisation functor $N$ is defined by
$$
N^{-1}(V)_n = \bigoplus_{\begin{smallmatrix} m+s=n\\ 0 \le \nu_1< \nu_2 < \ldots < \nu_s < n \end{smallmatrix}} \sigma_{\nu_s}\ldots \sigma_{\nu_2}\sigma_{\nu_1} V_m,
$$
where $\sigma_{\nu_s}\ldots \sigma_{\nu_2}\sigma_{\nu_1} V_m$ is just a copy of $V_m$, and we define $\sigma$ and $\pd$ on $N^{-1}(V)$ using the simplicial identities, subject to 
$$
\pd_i|_V = \left\{ \begin{matrix} d &i=0 \\ 0 &i>0. \end{matrix} \right.
$$
 Note that $N(N^{-1}(V)) \cong V$, and,  given $v \in V_n$, we write $N^{-1}v$ for 
$$
v \in (\id)V_n \subset N^{-1}(V)_n.
$$
\end{definition}

\begin{definition}
Given a local  complete pro-Artinian DG k-algebra $\cT_{\bullet}$, we can follow Quillen (\cite{QRat} Proposition 4.4) in defining the left adjoint to $N$, on complexes of algebras (adapted slightly for local algebras), by
$$
N^*(\cT)_n = F(N^{-1}(\m_{\cT}))_n/I_n,
$$
recalling that $F$ is our free functor for local pro-Artinian algebras, where $I_n$ is the ideal generated by
$$
\{N^{-1}(x) \nabla N^{-1}(y) - N^{-1}(x \wedge y)\colon\, x \in (\m_{\cT})_r, \, y \in (\m_{\cT})_s,\,r+s=n\},
$$
where $\wedge$ is the product on $\cT_{\bullet}$, and $\nabla$ is the shuffle product on $F(N^{-1}(\m_{\cT}))$, defined in terms of the (free) product on $F(N^{-1}(\m_{\cT}))$. We regard $(k \subset N^*(\cT)_n)$ as $\sigma_0^n(k \subset \cT_0)$ so that there is a copy of $N^{-1}(\cT)$ (and not just $N^{-1}(\m_{\cT})$) in $N^*(\cT)$.
\end{definition}

\begin{theorem}\label{qrat1}
If $\cT_{\bullet}$ is smooth, then so is $N^*(\cT)_*$.
\begin{proof}
This is essentially \cite{QRat} Proposition 4.4.
\end{proof}
\end{theorem}

\begin{theorem}\label{qrat2}
If $\cT_{\bullet}$ is smooth, then the unit $\eta:\cT_{\bullet} \to N(N^*(\cT))_{\bullet}$ is such that  the $H_i(\eta)$ are isomorphisms for all $i$. Thus $ \H_i(\cT_{\bullet}) \cong \pi_i(N^*(\cT)_*)$.
\begin{proof}
Essentially \cite{QRat} (Theorem 4.6 and the remark following). Note that the assumption of characteristic zero is vital here.
\end{proof}
\end{theorem}

It therefore seems that, if $\cT_{\bullet}$ is the smooth DG Hopf algebra corresponding to the DGLA $L^{\bullet}$, then a good candidate for the corresponding SDC $\cE(L)^*$ should be that represented by $N^*(\cT)$. However, we need first to define  Hopf maps on $N^*(\cT)$. Since the Hopf maps did not behave under $N$, it is unsurprising that they do not behave under $N^*$. We get round this by making use of the isomorphism $\bar{N} \cong N$.

\begin{lemma}
If $\cQ_*$ ``represents'' an SDC $E^*$, then, under the isomorphism
$$
N(\cQ)_{\bullet} \to \bar{N}(\cQ)_{\bullet},
$$
the map $\rho$ on $\bar{N}$ corresponds to $\Delta$ on $N$, where
$$
\Delta_{m,n}=((\id + \sum_{i=0}^{m-1}(-1)^{m-i}\sigma_i\pd_m)\ten \id) \circ \rho_{m,n}:
$$
$$
\begin{CD}
N(\cQ)_{\bullet} @>\Delta>> N(\cQ)_{\bullet}\hat{\ten} N(\cQ)_{\bullet}\\
@V\cong VV			@VV\cong V\\
\bar{N}(\cQ)_{\bullet}@>\rho>> \bar{N}(\cQ)_{\bullet}\hat{\ten}\bar{N}(\cQ)_{\bullet}.
\end{CD}
$$
On $N$, we may also recover $\rho$ from $\Delta$.
\begin{proof}
Observe that 
$$
\Delta(x) - \rho(x) \in D(\cQ)\hat{\ten}\cQ,
$$
 and verify that $\Delta(x) \in N(\cQ)\hat{\ten} N(\cQ)$.

Finally, note that we can apply the identities from Lemma \ref{*rules} to get
$$
\Delta_{m,n} = \rho_{m,n} +\sum_{i=0}^{m-1}(-1)^{m-i} (\sigma_i\ten \pd_0) \circ \rho_{m-1,n+1}.
$$
Hence
\begin{eqnarray*}
&&\Delta_{m,n} - \sum_{i=0}^{m-1}(-1)^{m-i} (\sigma_i\ten \pd_0) \circ \Delta_{m-1,n+1}\\
&=&  \rho_{m,n} +\sum_{i=0}^{m-1}(-1)^{m-i} (\sigma_i\ten \pd_0) \circ \rho_{m-1,n+1}
- \sum_{i=0}^{m-1}(-1)^{m-i} (\sigma_i\ten \pd_0)\circ \rho_{m-1,n+1}\\
&& -\sum_{i=0}^{m-1}\sum_{j=0}^{m-2}(-1)^{m-i}(-1)^{m-1-j}(\sigma_i\ten \pd_0)\circ (\sigma_j\ten \pd_0)\circ \rho_{m-2,n+2}\\
&=& \rho_{m,n} -\sum_{i=0}^{m-1}\sum_{j=0}^{m-2}(-1)^{m-i}(-1)^{m-1-j}(\sigma_i\sigma_j \ten  \pd_0\pd_0)\circ \rho_{m-2,n+2}\\
&=& \rho_{m,n}
\end{eqnarray*}
on $N$, since $\pd_0\pd_0=\pd_0\pd_1$, and $\id \ten \pd_1$ is zero on $\rho(N)$.
\end{proof}
\end{lemma}

\begin{definition}
Therefore, given a DGLA $L$, ``represented'' by the  smooth DG Hopf algebra $\cT_{\bullet}$ with comultiplication $\Delta$, we define  $\bar{N}^*(\cT)$  to be the algebra $N^*(\cT)$ (noting that this is smooth by Theorem \ref{qrat1}), with  Hopf maps defined on 
$$
\{N^{-1}(t):t \in \m_{\cT}\}\subset N^{-1}(\m_{\cT})
$$ 
by
$$
\rho_{m,n}:=\Delta_{m,n} - \sum_{i=0}^{m-1}(-1)^{m-i} (\sigma_i\ten \pd_0) \circ \Delta_{m-1,n+1},
$$
 extending to the rest of $N^{-1}(\m_{\cT})$ by the rules in Lemma \ref{*rules}, and extending to $N^*(\cT)$ by multiplicativity. Note that Lemma \ref{*shuffle} ensures that this is well-defined. Call the resulting SDC $\cE(L)^*$.
\end{definition}

\begin{theorem}\label{ltog} There are  canonical isomorphisms $ \H^i(L) \to\pi^i(\mathpzc{e}(L)^*) $, where $\mathpzc{e}(L)^*$ is the tangent space of $\cE(L)^*$, and  an equivalence $\Def_L \to \Def_{\cE(L)}$.
\begin{proof}
$ $
\begin{enumerate}
\item
$\mathpzc{e}(L)^* = \Der(\bar{N}^*(\cT)_*,k)$. From the definition of $\bar{N}^*(\cT)$, it follows that 
$$
\Der(\bar{N}^*(\cT)_*,k) \cong \{ \alpha \in \Hom(N^{-1}(\cT)^*,k) \colon \alpha(1)=0,\,\alpha(N^{-1}(x \wedge y))=0 \, \forall x,y \in \m_{\cT}\}.
$$
Now, one possible definition of cohomotopy groups is as the cohomology of the complex 
$$
\bar{N}^n(\mathpzc{e})=\mathpzc{e}^n/(\sum_{i=1}^n \pd^i \mathpzc{e}^{n-1}).
$$
Therefore
\begin{eqnarray*}
\pi^i(\mathpzc{e}(L)^*)&=&\H^i(\bar{N}(\mathpzc{e}(L))^{\bullet})\\
&=& \H^i(\bar{N}(\Der(\bar{N}^*(\cT)_*,k)))\\
&=& \H^i(\bar{N}(\{ \alpha \in \Hom(N^{-1}(\m_{\cT})_*,k) \colon \alpha(N^{-1}(x \wedge y))=0 \quad \forall x,y \in \m_{\cT}\}))\\
&=& \H^i(\{ \alpha \in \Hom(NN^{-1}(\m_{\cT})_{\bullet},k)\colon \alpha(N^{-1}(x \wedge y))=0 \quad \forall x,y \in \m_{\cT}\})\\
&=& \H^i(\Der(\cT_{\bullet},k))\\
&=&\H^i(L),
\end{eqnarray*}
since $NN^{-1}(\m_{\cT})_{\bullet}=(\m_{\cT})_{\bullet}$.

\item
First observe that $\bar{N}^*(\cT)_0= \cT_0$, with the same product. Thus $\gl=\cE(L)^0$. Call the identity of both these groups $e$.
Given $\alpha \in \mcl(A)$,
let 
$$
\omega: \bar{N}^*(\cT)_1 \to A
$$
 be the multiplicative map defined on the generators 
$$
N^{-1}(\cT)_1=\cT_1 \oplus \sigma_0 \cT_0 \quad\text{ by }\quad\omega = (\alpha ,0) + \omega_0,
$$
 where $\omega_0 = e \circ \pd_0$ [effectively this gives us $\omega = (\alpha ,e \circ (\sigma_0)^{-1})$]. The proof of Lemma \ref{gtol}, and the discussion comparing $N$ and $\bar{N}$, show that this is well-defined and satisfies the Maurer-Cartan equation. The gauge action of $\gl$ corresponds to the adjoint action of $E^0$, combining to give us a morphism
$$
\Def_L \to \Def_{\cE(L)}.
$$
That this is an equivalence follows from the first part, and Theorems \ref{SSC} to \ref{Man3}.
\end{enumerate}
\end{proof}
\end{theorem}

\begin{lemma}\label{gtoglg}
There is a quasi-isomorphism $E^* \to \cE(\cL(E))^*$ of SDCs, for any SDC $E^*$.
\begin{proof}
To see that there is a morphism $E^* \to \cE(\cL(E))^*$ is immediate. It is the map
$$
 \bar{N}^*\bar{N}(\cQ)_*  \to \cQ_*
$$
 corresponding to the canonical map $N^*N(\cQ_*) \to \cQ_*$. Moreover
$$
\pi_*(\bar{N}^*\bar{N}(\cQ)_*) \cong \H_*(\bar{N}(\cQ)_{\bullet})  \cong \pi_*(\cQ_*),
$$
the first isomorphism being a consequence of Theorem \ref{qrat2}.
Thus the map $\bar{N}^*\Sm\bar{N}(\cQ)_* \to \cQ_*$ induces an isomorphism on homotopy groups. Now, we may use the completeness of these simplicial algebras to show that they are free in the sense of \cite{QRat} Proposition I 4.4. Since free simplicial complexes are the analogue of projective resolutions, and we have a pair of homotopic free complexes, it follows that
$$
\Der(\cQ_*,k) \to \Der(\bar{N}^*\bar{N}(\cQ)_*,k)
$$
gives an isomorphism on cohomotopy groups, i.e.
$$
\pi^i(\mathfrak{e}^*) \cong \pi^i(\mathpzc{e}(\cL(E))^*),
$$
the proof being the same as the result in \cite{Q} that the quasi-isomorphism class of the tangent complex is independent of the choice of cofibrant resolution.
Hence $E^* \to \cE(\cL(E))^*$ is a quasi-isomorphism of SDCs.
\end{proof}
\end{lemma}

\begin{corollary}\label{ltolgl}
For any DGLA $L^{\bullet}$, there is a quasi-isomorphism $L^{\bullet} \to \cL(\cE(L))^{\bullet}$.
\begin{proof}
By Lemma \ref{gtoglg}, there is a quasi-isomorphism of SDCs
$$
\cE(L)^* \to \cE(\cL(\cE(L)))^*.
$$
Now, Theorem \ref{ltog} gives a canonical isomorphism $\bar{N}(\mathpzc{e}(L))^{\bullet} \cong L^{\bullet}$, so we may regard $\bar{N}(\mathpzc{e}(L))^{\bullet}$ as a DGLA, with bracket 
$$
[\alpha,\beta]= (\alpha \ten \beta)\circ \Delta \mp (\beta\ten \alpha)\circ \Delta.
$$ 
(Note that $\bar{N}(\mathfrak{e}^*)$, for the tangent space $\mathfrak{e}^*$ of an arbitrary SDC $E^*$ will \emph{not} be closed under this bracket.) 

The quasi-isomorphism of SDCs then gives us a quasi-isomorphism of DGLAs
$$
L^{\bullet} \cong \bar{N}(\mathpzc{e}(L))^{\bullet} \xra{\sim} \bar{N}(\mathpzc{e}(\cL(\cE(L))))^{\bullet} \cong \cL(\cE(L))^{\bullet}.
$$
\end{proof}
\end{corollary}

\begin{theorem}\label{gtol2}
 The maps $\pi^i(\mathfrak{e}^*) \to \H^i(\cL(E)^{\bullet})$ in Lemma \ref{hgtol} are isomorphisms, and hence the morphism $\Def_E \to \Def_{\cL(E)}$ in Lemma \ref{gtol} is an equivalence.
\begin{proof}
We can factor the isomorphism 
$$
\pi^i(\mathfrak{e}^*) \xra{\cong} \pi^i(\mathpzc{e}(\cL(E))^*) 
$$
from Lemma \ref{gtoglg} as 
$$
\pi^i(\mathfrak{e}^*) \to \H^i(\cL(E)^{\bullet}) \xra{\cong} \pi^i(\mathpzc{e}(\cL(E))^*),
$$
the second isomorphism coming from Theorem \ref{ltog}. Thus 
$$
\pi^i(\mathfrak{e}^*) \to \H^i(\cL(E)^{\bullet})
$$ 
must be an isomorphism, and by Theorems \ref{SSC} to \ref{Man3} and Lemma \ref{gtol},
$$
\Def_E \to \Def_{\cL(E)}
$$
must be an equivalence.
\end{proof}
\end{theorem}

\subsubsection{An explicit description}\label{explicit}

Given a DGLA $L$, first form the denormalised cosimplicial complex
$$
D^nL:= (\bigoplus_{\begin{smallmatrix} m+s=n \\ 0 \le j_1 < \ldots < j_s \le n \end{smallmatrix}} \pd^{j_s}\ldots\pd^{j_1}L^m)/(dv \sim \sum_{i=0}^{n+1}(-1)^i \pd^i v \quad \forall v \in L^n),
$$
where    we define the $\pd^j$ and $\sigma^i$ using the simplicial identities, subject to the condition that $\sigma^i L =0$.

We now have to define the Lie bracket $\llbracket,\rrbracket$ from $D^nL \ten D^nL$ to $D^n L$. On the conormalised complex $N^n(DL\ten DL) \le D^nL \ten D^nL$, we define this as the composition
$$
N^n(DL\ten DL) \xra{\sum \nabla^{pq}} \bigoplus_{p+q=n}L^p\ten L^q \xra{[,]} L^n, 
$$
where 
$$
\nabla^{pq}(v\ten w)=\sum_{ (\mu, \nu) \in \mathrm{Sh}(p,q)}(-1)^{(\mu,\nu)} (\sigma^{\nu_1}\ldots \sigma^{\nu_q}v)\ten (\sigma^{\mu_1}\ldots \sigma^{\mu_p}w).
$$

We then  extend this bracket to the whole of $D^nL\ten D^n L$ by setting
$$
\llbracket \pd^ia,\pd^ib \rrbracket:=\pd^i \llbracket a,b \rrbracket.
$$

Now set
$$
\cE(L)^n(A)= \exp(D^n(L)\ten \m_A),
$$
making $\cE(L)$ into a cosimplicial complex of group-valued functors. To make it an SDC, we must define a $*$ product. We do this as the Alexander-Whitney cup product
$$
g*h = (\pd^{m+n}\ldots \pd^{m+2}\pd^{m+1}g)\cdot (\pd^0)^m h,
$$
 for $g \in \cE(L)^m,\,h \in \cE(L)^n$.


\newpage
\bibliographystyle{alpha}
\addcontentsline{toc}{section}{Bibliography}
\bibliography{references.bib}
\end{document}